\documentclass[10pt]{article}
\usepackage{graphicx}
\usepackage{amsfonts}
\usepackage{amsmath}
\usepackage{amssymb}
\usepackage{eufrak}
\usepackage{color}
\usepackage{algorithm}
\usepackage{algorithmic}
\usepackage{rotating}
\usepackage{textcomp}
\usepackage{hyperref}
\def\Prox{{\hbox{\rm Prox}}}
\def\cS{{\mathcal{S}}}
\def\Tr{{\hbox{\rm Tr}}}
\def\cA{{\mathcal{A}}}
\def\cL{{\cal L}}
\def\Ent{{\hbox{\rm Ent}}}
\def\cH{{\cal H}}
\def\bE{{\mathbb{E}}}
\def\cN{{\cal N}}
\def\cM{{\cal M}}

\def\Opt{{\hbox{\rm Opt}}}
\def\Prob{{\hbox{\rm Prob}}}
\def\cC{{\cal C}}
\newtheorem{thm}{Theorem}[section]
\newtheorem{cor}{Corollary}[section]
\newtheorem{lem}[cor]{Lemma}
\newtheorem{prop}[cor]{Proposition}

\newtheorem{rem}{Remark}[section]
\newtheorem{ass}{Assumption}[section]

\newtheorem{myalgorithm}{Algorithm}

\newcommand{\R}{\mathbb R}

\newcommand{\proof}{

\vspace{-2ex}

\noindent\textbf{\textit{Proof}}\\}

\newcommand{\qed}{${}$\hfill \rule{2mm}{2mm}}

\begin{document}
\title{A randomized Mirror-Prox method for solving \\structured large-scale matrix saddle-point problems}
\author{Michel Baes\footnote{Institute for Operations Research, ETH Z\"urich, R\"amistrasse 101, 8092 Z\"urich, Switzerland,\newline michel.baes@ifor.math.ethz.ch.} , Michael B\"urgisser\footnote{Institute for Operations Research, ETH Z\"urich, R\"amistrasse 101, 8092 Z\"urich, Switzerland,\newline michael.buergisser@ifor.math.ethz.ch.} , Arkadi Nemirovski\footnote{Georgia Institute of Technology, Atlanta, Georgia 30332, USA, nemirovs@isye.gatech.edu.}}

\maketitle

\begin{abstract}
In this paper, we derive a randomized version of the Mirror-Prox method for solving some structured matrix saddle-point problems, such as the maximal eigenvalue minimization problem. Deterministic first-order schemes, such as Nesterov's Smoothing Techniques or standard Mirror-Prox methods, require the exact computation of a matrix exponential at every iteration, limiting the size of the problems they can solve. Our method allows us to use stochastic approximations of matrix exponentials. We prove that our randomized scheme decreases significantly the complexity of its deterministic counterpart for large-scale matrix saddle-point problems. Numerical experiments illustrate and confirm our theoretical results.
\end{abstract}


\section{Introduction}
Large-scale semidefinite programming attracts substantial research
efforts nowadays. A vast set of applications can be modeled as such
optimization problems, and many strategies have been studied
theoretically and implemented in excellent softwares.

Arguably, general purpose semidefinite methods suffer from an intrinsic
drawback. They forbid themselves, for the sake of generality, to
exploit explicitly some structural features of the particular
instance they are given to solve, hampering the resolution of very
large-scale problems.

As a result, we are witnessing the development of special-purpose
algorithms, designed for particular subclasses of semidefinite optimization problems, where the
utilization of their specific structure is instrumental for aiming at
large size problems. In this paper, we are addressing the problem of
minimizing the maximal eigenvalue of an affine combination of given
symmetric matrices, plus a linear function of the coefficients of
this affine combination. Different strategies have been devised to
deal specifically with the maximal eigenvalue minimization problem.
Among the first investigated techniques, Bundle methods were
introduced in \cite{Helmberg_bundle,oustry}, and subsequently refined in a
number of further papers. Theoretical results on Bundle methods
concern mainly asymptotic convergence properties; to the best of our
knowledge, no complexity guarantees have been obtained so far in this
direction.

Another strategy has been discovered in \cite{nesterov:coreDP73/2004}
when Nesterov showed how his Smoothing Techniques can be specialized
to the maximal eigenvalue minimization problem. As a result, he obtained
worst-case complexity guarantees for his method: if $\epsilon>0$ is
the desired absolute accuracy on the objective value, $A_1,\ldots,A_m$ are real symmetric $n\times n$-matrices, $\Delta_m\subseteq\R^m$ is the $(m-1)$-dimensional simplex, and $\lambda_{\max}(A)$ denotes the maximal eigenvalue of any real symmetric matrix $A$, his algorithm
solves the problem
\begin{eqnarray}
 \label{eq:problem_nesterov}
  \min \left\{\lambda_{\max}\left(\sum_{j=1}^mA_jx_j\right):\ x\in\Delta_m\right\}
\end{eqnarray}
in
$\mathcal{O}((n^3+mn^2)\max_{j}\lambda_{\max}(\left|A_j\right|)\sqrt{\ln (n)\ln (m)}/\epsilon)$ elementary
operations, where $\left|A\right|=\sqrt{A^2}$ for any real symmetric matrix $A$. Almost simultaneously, the paper
\cite{Nemirovski_Prox_method_saddle} develops a Mirror-Prox method,
which can be particularized as well for our problem, and obtains
equivalent complexity results.  Two papers by Warmuth et al.
\cite{Tsuda:matrix_exponentiated,Warmuth:online_variance} present a scheme - called the Matrix Exponentiated Gradient Update method - that can basically be applied to problem (\ref{eq:problem_nesterov}). A form of this algorithm was independently discovered by Arora and Kale \cite{Arora_Kale_07}. The methods of Arora et al. and Warmuth et al. essentially reduce to a subgradient method (provided that we adapt them to our problem),
but with a worse complexity guarantee: in order to find a solution to problem (\ref{eq:problem_nesterov}) with absolute accuracy $\epsilon$, these methods need
$\mathcal{O}(\max_j\lambda_{\max}^2(\left|A_j\right|)\ln (n)/\epsilon^2)$ iterations. Each of these iterations requires the computation of a matrix exponential and further operations with a cost not exceeding $\mathcal{O}(mn^2)$.

In a nutshell, the methods introduced in \cite{Arora_Kale_07,Tsuda:matrix_exponentiated,Warmuth:online_variance} present the same computational bottleneck as Smoothing Techniques and the Mirror-Prox method when
applied to our problem: at every iteration, all these schemes require the
determination of a symmetric matrix's exponential. Several efforts
have been carried out to reduce the iteration computation cost. In \cite{d'Aspremont:Approx_gradients}, d'Aspremont
analyzes the possibility of using approximate gradients, and thereby
approximate matrix exponentials in Nesterov's Smoothing Techniques. In
\cite{Juditsky08_solvingvariational}, Mirror-Prox methods for variational
problems where extended to situations where only some stochastic
information is available from the instance to solve. These methods
were particularized to the maximal eigenvalue minimization problem
where all the input matrices share the same block-diagonal pattern.
Albeit the problem is completely deterministic, an artificial
randomization was introduced in the oracle of the method, which
reduced the iteration cost while retaining some probabilistic
guarantees on the output of the algorithm. Finally, Arora and Kale \cite{Arora_Kale_07} obtain, by approximating the rows of $\exp(X/2)$, a substitute for the exact Gram matrix $\exp(X)$, where $X$ is some real symmetric $n\times n$-matrix. The rows of $\exp(X/2)$ are approximated by projecting an appropriate truncation of the exponential Taylor series approximation on, roughly speaking, $\mathcal{O}(1/\epsilon^2)$ random directions.

In this paper, we apply the general results of
\cite{Juditsky08_solvingvariational} to analyze another randomization
strategy for computing matrix exponentials, which is also based on a vector sampling
and on an appropriate truncation of the exponential Taylor series. Whereas we consider the same number of terms in the Taylor series approximation of the matrix exponential as Arora and Kale \cite{Arora_Kale_07} do, we can significantly reduce the number of required random vectors: roughly speaking, we project the truncated Taylor series on $\mathcal{O}(1/\epsilon)$ random directions.

The approximation strategy developed in this paper proves to be theoretically efficient for large-scale problems. In theory, it outperforms all its competitors on a reasonably large set of instances, described by the size of the input matrices, their number, their sparsity, and the desired accuracy. Our
theoretical conclusions are demonstrated by numerical evidence: for solving problems (up to a relative accuracy of $0.2$\%) that involve a hundred matrices of size $800\times 800$, the Mirror-Prox method equipped with our randomization procedure requires on average about, roughly speaking, half of the CPU time needed by the Mirror-Prox method with exact computations.

The paper is organized as follows. Section 2 contains the necessary
notational conventions and a brief recall on existing results on
Mirror-Prox methods for general convex problems with approximate
oracle. We particularize these considerations in Section 3 to slightly structured matrix
saddle-point problems and we analyze the stochastic exponential
approximation strategy briefly described above for computing an
approximate oracle. In Section 4, we derive the complexity of solving
the maximal eigenvalue minimization problem. In Section 5, we test
our method for solving large-scale eigenvalue optimization problems,
comparing its efficiency with the provably best purely deterministic
method in terms of worst-case complexity.

\section{Mirror-Prox methods with approximate first-order information: a review}

Let $E$ be a Euclidean space with inner product $\left\langle  \cdot,\cdot\right\rangle$. We endow $E$ with a norm $\left\|\cdot\right\|$, which may differ from the one that is induced by this inner product. The {\sl conjugate norm} $\left\|\cdot\right\|_{\ast}$ to $\left\|\cdot\right\|$ is defined as:
\[\left\|w\right\|_\ast:=\max_{x\in E}\left[\left\langle w,x\right\rangle:\ \left\|x\right\|=1\right].\]

\subsection{Variational inequalities and saddle-point problems}

\paragraph{Variational inequalities.} Let $Q$ be a non-empty convex compact subset of $E$, and let $F:Q\rightarrow E$ be a Lipschitz continuous monotone mapping with Lipschitz constant $L>0$:
\begin{eqnarray*}
 \label{ass:monotone_Lipschitz}
\left\|F(z)-F(z')\right\|_\ast\leq L\left\|z-z'\right\| &\qquad& \forall\ z,z'\in Q\cr
\left\langle F(z)-F(z'),z-z'\right\rangle\geq 0 &\qquad& \forall\ z,z'\in Q.
\end{eqnarray*}
The variational inequality associated with the set $Q$ and the operator $F$ reads as follows:
\begin{equation}
\label{eq:vi}
\hbox{find $z^\ast\in Q$ such that $\left\langle F(z),z^\ast-z \right\rangle \leq 0$ for all $z\in Q$.}
\end{equation}
In the sequel, in order to measure the inaccuracy of a point $\bar{z}\in Q$ as a candidate solution to (\ref{eq:vi}), we use the {\sl dual gap function}
\[\epsilon(\bar{z}):=\max_{z\in Q}\left\langle F(z),\bar{z}-z \right\rangle. \]
For $z\in Q$, we clearly have $\epsilon(z)\geq0$, and $\epsilon(z)=0$ if and only if $z$ solves (\ref{eq:vi}).
\paragraph{Saddle point problems.} Assume that $E:=E_1\times E_2$ for Euclidean spaces $E_1$ and $E_2$, and that $Q:=Q_1\times Q_2$ is non-empty with two convex compact sets $Q_1\subset E_1$ and $Q_2\subset E_2$. Let $\phi:Q_1\times Q_2\rightarrow\mathbb{R}$ be a convex-concave function. We restrict ourselves to functions $\phi$ that are differentiable with Lipschitz continuous gradient. The function $\phi(\cdot,\cdot)$ is associated with the saddle-point problem
\begin{eqnarray}
\label{eq:saddle-point}
\min_{x\in Q_1}\max_{y\in Q_2} \phi(x,y).
\end{eqnarray}
Due to the standard Minimax Theorem in Convex Analysis (see Corollary 37.3.2 in \cite{pm:Rockafellar:70}), we have the following pair of primal-dual convex optimization problems:
\[\min_{x\in Q_1}\left[\overline{\phi}(x):=\max_{y\in Q_2}\phi (x,y)\right] =\max_{y \in Q_2}\left[\underline{\phi}(y):=\min_{x\in Q_1}\phi (x,y)\right].\]
It is well known that the solutions to the saddle point problem (\ref{eq:saddle-point}) are exactly the pairs $(x_*,y_*)$ comprised of optimal solutions to the above two optimization problems, and that these pairs are exactly the solutions to the variational inequality given by $Q=Q_1\times Q_2$ and the monotone operator \[F(x,y)=\left(\frac{\partial \phi(x,y)}{\partial x}; -\frac{\partial \phi(x,y)}{\partial y}\right):Q\to E_1\times E_2.\] We quantify the accuracy of a candidate solution $\bar{z}=(\bar{x},\bar{y} )\in Q$ to the saddle point problem (\ref                                             {eq:saddle-point}) by the value of the corresponding duality gap
\[\epsilon^{sad}(\bar{z}):=\overline{\phi}(\bar{x})-\underline{\phi}(\bar{y})=\max_{v\in Q_2}\phi(\bar{x},v)-\min_{u\in Q_1}\phi(u,\bar{y}).\]

Due to the convex-concave structure of $\phi$, any point $\bar{z}=\left(\bar{x},\bar{y}\right)\in Q_1\times Q_2$ constitutes an $\epsilon^\text{sad}(\bar{z})$-approximate solution to the variational inequality that is associated with $Q=Q_1\times Q_2$ and the above $F$:
\begin{eqnarray*}
\epsilon^{sad}(\bar{x},\bar{y}) &=&  \max_{(x,y)\in Q_1\times Q_2}\left[\phi(\bar{x},y)-\phi(x,y)+\phi(x,y)-\phi(x,\bar{y})\right]\cr
 &\geq&\max_{(x,y)\in Q_1\times Q_2}\left[\left \langle {\partial\phi(x,y)\over\partial x},\bar{x}-x\right\rangle+\left\langle -{\partial\phi(x,y)\over\partial y},\bar{y}-y\right\rangle\right]=\max_{z\in Q} \left\langle F(z),\bar{z}-z\right\rangle=\epsilon(\bar{x},\bar{y}).
\end{eqnarray*}

\subsection{Mirror-Prox algorithm: preliminaries}\label{sectMP:preliminaries}
In its basic form, the Mirror Prox (MP) algorithm is aimed at solving variational inequalities on a convex compact subset $Q$ of a Euclidean space $E$ equipped with a norm $\|\cdot\|$. The setup for the algorithm is given by a {\sl distance-generating function} (d.-g.f.) $\omega$: $Q\to\R$ which possesses the following properties:
\begin{itemize}
 \item $\omega$ is continuous and convex on $Q$. In particular, the domain $Q^o=\{x\in Q: \partial \omega(x)\neq\emptyset\}$ of the subdifferential of $\omega$ is nonempty.
 \item $\omega$ is regular on $Q^o$, i.e., the subdifferential $\partial\omega(\cdot)$ admits a {\sl continuous} selection  $\omega'(\cdot)$ on $Q^o$.
\item The function $\omega$ is strongly convex modulus 1 with respect to $\|\cdot\|$:
\[\left\langle \omega'(z) -\omega'(y),z-y\right\rangle \geq \left\|z-y\right\|^2\qquad \forall\ y,z\in Q^o.\]
In the sequel, we refer to the latter property as the {\sl compatibility} of the d.-g.f. $\omega(\cdot)$ and the norm $\|\cdot\|$.
\end{itemize}
Furthermore, we suppose that we choose $\omega$ such that we can easily solve problems of the form:
\begin{eqnarray}
\label{eq:easy_solvable}
\min_{z\in Q}\left[ \omega(z)-\left\langle e,z\right\rangle  \right],\qquad e\in E.\end{eqnarray}
\begin{rem}\label{rem:z_e} The optimal solution $z_e$ to (\ref{eq:easy_solvable}) clearly exists, is unique by continuity and strong convexity of $\omega$, and belongs to $Q^o$ (indeed, by evident reasons $e\in\partial\omega(z_e)$). From regularity of $\omega(\cdot)$, it immediately follows that
\begin{equation}\label{eqArkadi}
\langle \omega'(z_e)-e,z-z_e\rangle \geq0\qquad \forall\ z\in Q.
\end{equation}
\end{rem}
A d.-g.f. $\omega(\cdot)$ gives rise to several entities:
\begin{itemize}
\item  The {\sl $\omega$-center} $z^\omega:=\arg\min_{z\in Q} \omega(z)$ of $Q$.
\item The {\sl Bregman distance} $V_z(w)=\omega(w)-\omega(z)-\langle\omega'(z),w-z\rangle$, where $z\in Q^o$ and $w\in Q$. By strong convexity, we have:
\begin{equation}\label{eqArkadi1}
 V_z(w)\geq {1\over 2}\|w-z\|^2\qquad \forall (z\in Q^o,w\in Q).
\end{equation}
\item The {\sl $\omega$-diameter $\Omega$ of $Q$}, which is defined as:
\[
\Omega:=\sqrt{2\max_{z\in Q} V_{z^\omega}(z)}\leq\sqrt{2\left(\max_{z\in Q}\omega(z)  -\min_{z\in Q}\omega(z)\right)},
\]
where the concluding inequality follows from the fact that $V_{z^\omega}(z)\leq\omega(z)-\omega(z^\omega)$ due to $\langle\omega'(z^\omega),z-z^\omega\rangle\geq0$ for every $z\in Q$, see (\ref{eqArkadi}). Further, by (\ref{eqArkadi1}), we have:
\begin{equation}\label{eqArkadi2}
 \|w-z^\omega\|\leq\Omega\text{ for any }w\in Q,\hbox{\ whence\ } D:=\max\limits_{w,z\in Q}\|w-z\|\leq 2\Omega.
\end{equation}
\item For parameter $z\in Q^o$, we define the {\sl Prox-mapping} as:
$$
\Prox_z(\xi)=\arg\min_{w\in Q} \left[V_z(w)+\langle\xi,w\rangle\right]: E\to Q^o
$$
(the $\arg\min$ in question indeed belongs to $Q^o$, see Remark \ref{rem:z_e}).
\end{itemize}

\subsection{Mirror-Prox algorithm with noisy first-order information}

The prototype MP algorithm \cite{Nemirovski_Prox_method_saddle} is aimed at solving variational inequality (\ref{eq:vi}) when exact values of $F$ are available. In this paper, we use a modification of the original MP scheme, the {\sl Stochastic Mirror-Prox} (SMP) algorithm proposed and investigated in \cite{Juditsky08_solvingvariational}, which operates with noisy estimates of $F$. Specifically, the algorithm has access to a {\sl Stochastic Oracle}: at the $t$-th call of the oracle, $z_t\in Q^o$ being the query point, the oracle returns an estimate $\hat{F}_{\xi_t}(z_t)$ of $F(z_t)$. Here, $\xi_t$ is the $t$-th realization of the ``oracle's noise'', which is modeled as a random vector $\xi$, and $\hat{F}_{\xi}(z)$ is a Borel function of $\xi$ and $z$. We assume that the realizations $\left(\xi_t\right)_{t\geq 1}$ of the random vector $\xi$ are independent. From now on, we set $\xi_{[t]}=(\xi_1,\xi_2,...,\xi_t)$. The algorithm is as follows.

\begin{myalgorithm}\label{generalMP}
\label{alg:mirror_prox}{\rm[Mirror-Prox method with noisy first-order information]\\
1: Choose the number of iterations $T$.
 Set $z_0=z^\omega\in Q^o$.\\
2: {\bf for} $1\leq t\leq T$ {\bf do}\\
3: Given $z_{t-1}\in Q^o$, choose (a deterministic) $\gamma_t>0$ such that:
\begin{eqnarray}
\label{eq:stepsize_bound}
\gamma_t\leq \frac{1}{\sqrt{2}L}.
\end{eqnarray}
4: Call Stochastic Oracle with query point $z_{t-1}$ and receive $\eta_t:=\hat{F}_{\xi_{2t-1}}(z_{t-1})$. \\
5: Set $w_t=\Prox_{z_{t-1}}(\gamma_t\eta_t)\in Q^o$.\\
6: Call Stochastic Oracle with query point $w_t$ and receive $\zeta_t:=\hat{F}_{\xi_{2t}}(w_t)$. \\
7: Set $z_t=\Prox_{z_{t-1}}(\gamma_t\zeta_t)\in Q^o$.\\
8: {\bf end for}\\
9: Return $z^T:=\left(\sum_{t=1}^T\gamma_t\right)^{-1}\sum_{t=1}^T\gamma_t w_t.$
}
\end{myalgorithm}
Note that $z_t$ is a deterministic function of $\xi_{[2t]}$, while $w_t$ is a deterministic function of $\xi_{[2t-1]}$. In order to show expected convergence of Algorithm \ref{alg:mirror_prox}, we need to define the following quantities:
\begin{equation}\label{martingale}
\begin{array}{rcl}
 \mu_{{z}_{t-1}}&:=&\mathbb{E}_{\xi_{2t-1}}\left\{\hat{F}_{\xi_{2t-1}}({z}_{t-1})|\ \xi_{[2t-2]}\right\}-F({z}_{t-1}),\cr
\mu_{{w}_t}&:=&\mathbb{E}_{\xi_{2t}}\left\{\hat{F}_{\xi_{2t}}({w}_{t})| \xi_{[2t-1]}\right\}-F({w}_t),\cr
\sigma_{{z}_{t-1}}&:=&\hat{F}_{\xi_{2t-1}}({z}_{t-1})-\mathbb{E}_{\xi_{2t-1}}\left\{\hat{F}_{\xi_{2t-1}}({z}_{t-1})|\ \xi_{[2t-2]}\right\},\cr
\sigma_{{w}_t}&:=&\hat{F}_{\xi_{2t}}({w}_{t})-\mathbb{E}_{\xi_{2t}}\left\{\hat{F}_{\xi_{2t}}({w}_{t})| \xi_{[2t-1]}\right\},
\end{array}
\end{equation}
where  $1\leq t\leq T$. Note that we define $\mathbb{E}_{\xi_{1}}\left\{\hat{F}_{\xi_{1}}({z}_{0})|\ \xi_{[0]}\right\}:=\mathbb{E}_{\xi_{1}}\left\{\hat{F}_{\xi_{1}}({z}_{0})\right\}$. Note also that $\sigma_{z_{t-1}}$ and $\sigma_{w_t}$ are martingale differences. The following result is proven in the Appendix.

\begin{thm}
\label{thm:convergence} Let
\begin{equation}\label{epsilonT}
\epsilon^T=\frac{\frac{\Omega^2}{2}+\Omega\| \sum_{t=1}^T\gamma_t\sigma_{{w}_t}\|_*+\sum_{t=1}^T \left\{\gamma_t D \|\mu_{{w}_t}\|_*+{2\gamma_t^2}\left(\|\sigma_{{w}_t}-\sigma_{{z}_{t-1}}\|_*^2+\|\mu_{{w}_t}-\mu_{{z}_{t-1}}\|_*^2\right)\right\}}
{\sum_{t=1}^T\gamma_t}.
\end{equation}
We have:
\[\mathbb{E}_{\xi_{[2T]}}\left\{\epsilon(z^T)\right\}\leq \mathbb{E}_{\xi_{[2T]}}\left\{\epsilon^T\right\}.\]
Moreover, for an operator $F$ that is associated with saddle-point problem (\ref{eq:saddle-point}), the following inequality holds:
\[\mathbb{E}_{\xi_{[2T]}}\left\{\epsilon^\text{sad}(z^T)\right\}\leq \mathbb{E}_{\xi_{[2T]}}\left\{\epsilon^T\right\}.\]
\end{thm}

\section{Mirror-Prox algorithm for matrix saddle-point problems}

\subsection{Matrix saddle-point problems}
The problem of primary interest in this paper is the Eigenvalue Minimization problem
\begin{equation}\label{eqArkadi4}
\Opt=\min\limits_{x\in Q_1}\left[\lambda_{\max}(\cA(x)+B)+\langle c,x\rangle\right],
\end{equation}
where:
\begin{itemize}
\item $Q_1$ is a convex compact subset of the space $E_1=\R^{m}$ equipped with the standard inner product $\langle x,y\rangle=x^Ty$;
\item $\cA(x)=\sum_{i=1}^{m}x_iA_i$ is a linear mapping from $\R^m$ into the space $E_2=\cS^n$ of symmetric $n\times n$ matrices (so that $A_1,...,A_m\in\cS^n$), and $B\in \cS^n$;
\item $\lambda_{\max}(A)$ stands for the maximal eigenvalue of a symmetric matrix $A$;
\item $c\in\mathbb{R}^m$.
\end{itemize}
We equip $E_2=\cS^n$ with the Frobenius inner product $\langle X,Y\rangle_F=\Tr(XY)$. Denoting by $\Delta^M_n$  the {\sl spectahedron} $\{Y\in \cS^n:Y\succeq0,\Tr(Y)=1\}$ and observing that $\lambda_{\max}(A)=\max\limits_{Y\in\Delta_n^M} \Tr(YA)$, we can reformulate (\ref{eqArkadi4}) as the following bilinear saddle point problem:
\begin{eqnarray}
\label{eq:matrix_saddle-point}\label{eqArkadi44}
\Opt:=\min_{x\in Q_1}\max_{Y\in \Delta_n^M}\phi (x,Y),\qquad \phi(x,Y)=\left\langle B,Y\right\rangle_F+\left\langle \cA(x),Y\right\rangle_F+\left\langle c,x\right\rangle.
\end{eqnarray}
The associated operator $F$ is:
\begin{equation}
 \label{def:operatorF}
F(x,Y):=\left[F_x(Y):=\frac{\partial \phi(x,Y)}{\partial x}=\cA^\ast(Y)+c;
F_Y(x):=-\frac{\partial \phi(x,Y)}{\partial Y}=-\cA(x)-B\right],
 \end{equation}
 where the linear mapping $\cA^\ast(Y)=[\Tr(A_1Y);\Tr(A_2Y);...;\Tr(A_mY)]:\cS^n\to \R^m$ is conjugate to the mapping $x\mapsto \cA(x)$. \par
 We are about to solve the Eigenvalue Minimization problem by applying to (\ref{eqArkadi44}) the Mirror-Prox algorithm. While the problem is fully deterministic, we intend to use an appropriately constructed Stochastic Oracle, which is computationally cheaper than the exact deterministic oracle; our ultimate goal is to demonstrate that the resulting SMP algorithm significantly outperforms its deterministic counterpart in a meaningful range of problem's sizes.
 \par
 We start with presenting the algorithm's setup.

\subsection{Algorithm's setup}\label{sect:setupMSP}
We assume that the space $E_1=\R^m$ is equipped with a norm $\|\cdot\|_{x}$, the conjugate norm being $\|\cdot\|_{x,*}$, and that $Q_1$ is equipped with a d.-g.f. $\omega_x(x)$ that is compatible with $\|\cdot\|_{x}$. We denote by $x^{\omega_x}$ and  $\Omega_x$ the $\omega_x$-center of $Q_1$ and the $\omega_x$-diameter of $Q_1$, respectively; see Section \ref{sectMP:preliminaries}.\par
We equip the space $E_2=\cS^n$ with the trace-norm $\|W\|_Y:=\sum_{i=1}^n|\lambda_i(W)|$, where the vector $\lambda(W)=[\lambda_1(W);\lambda_2(W);...;\lambda_n(W)]$ consists of the eigenvalues $\lambda_1(W)\geq\ldots\geq\lambda_n(W)$ of $W\in\cS^n$. As it is well-known, the conjugate norm is the usual spectral norm $\|W\|_{Y,*}=\max_{1\leq i\leq n}|\lambda_i(W)|$. Further, we equip the spectahedron $Q_2:=\Delta_n^M$ with the matrix entropy d.g.-f.:
\begin{equation}\label{Ent}
\omega_Y(Y)=\Ent(Y):=\sum_{i=1}^n\lambda_i(Y)\ln(\lambda_i(Y)).
\end{equation}
 As shown in \cite{nesterov:coreDP73/2004}, see also \cite{Ben-tal_nemirovski_level}, this d.-g.f. is compatible with $\|\cdot\|_{Y}$, and as it is immediately seen, the corresponding center and diameter of $Q_2$ are as follows:
\begin{equation}\label{parameters}
Y^{\omega_Y}={1\over n}I_n;\quad \Omega_Y=\sqrt{2\ln(n)}.
\end{equation}
Finally, we equip the embedding space $E=E_1\times E_2$ of the domain $Q=Q_1\times Q_2$ of (\ref{eqArkadi44}) with the norm
\begin{equation}\label{norm}
\|(x,Y)\|=\sqrt{{1\over \Omega_x^2}\|x\|_x^2+{1\over\Omega_Y^2}\|Y\|_Y^2},
\end{equation}
implying that the conjugate norm is
\begin{equation}\label{conjnorm}
\|(x,Y)\|_*=\sqrt{\Omega_x^2\|x\|_{x,*}^2+\Omega_Y^2\|Y\|_{Y,*}^2}.
\end{equation}
The domain $Q=Q_1\times Q_2$ of (\ref{eqArkadi44}) is equipped with the d.-g.f.
\begin{equation}\label{dgf}
\omega(x,Y)={1\over \Omega_x^2}\omega_x(x)+{1\over \Omega_Y^2}\omega_Y(Y);
\end{equation}
it is immediately seen that this indeed is a d.-g.f. for $Q$ compatible with $\|\cdot\|$, and that the corresponding diameter of $Q$ is $\Omega=\sqrt{2}$, while the $\omega$-center of $Q$ is $z^\omega=(x^{\omega_x},Y^{\omega_Y})$.
\par
Finally, let $\cL$ be (an upper bound on) the norm of the linear mapping $x\mapsto\cA(x):E_1\to E_2$ induced by the norms $\|\cdot\|_x$ and $\|\cdot\|_{Y,*}$ on the argument and the image spaces:
\begin{equation}\label{eqArkadi66}
\forall x\in E_1: \|\cA(x)\|_{Y,*}\leq \cL\|x\|_x.
\end{equation}
It is immediately seen that the affine monotone operator $F$ associated with (\ref{eqArkadi44}) (see (\ref{def:operatorF})) satisfies:
\begin{equation}\label{eqArkadi6}
\forall (z,z'\in E=E_1\times E_2): \|F(z)-F(z')\|_*\leq L\|z-z'\|,\qquad\text{where }L:=\Omega_x\Omega_Y\cL.
\end{equation}

\subsection{Randomized Mirror-Prox method for (\ref{eq:matrix_saddle-point})}

\subsubsection{Randomization: motivation and strategy} \label{subs:strategy}

With the outlined setup,  when applying the deterministic MP algorithm (i.e., Algorithm \ref{alg:mirror_prox} with precise information:  $\hat{F}_{\xi_t}\equiv F$) to the variational inequality associated with the saddle point problem (\ref{eqArkadi44}), the computational effort at iteration $t$ is dominated by the necessity

\begin{itemize}
\item[(A)] to compute the value of $F$ at two points, namely at the points $\bar{z}=z_{t-1}=(\bar{x},\bar{Y})\in Q^o$ and $\bar{w}=w_t\in Q^o$;
\end{itemize}
and
\begin{itemize}
\item[(B)] to compute the value of the prox-mapping $\Prox_{\bar{z}}(\zeta)=\arg\min_{w\in Q}\left[\omega(w)+\langle \zeta-\omega'(\bar{z}),w-\bar{z}\rangle\right]$ at two different points $\zeta\in E$.
\end{itemize}
With our d.-g.f. $\omega$ that is ``separable'' with respect to the $x$- and to the $Y$-component of $\bar{z}$, task (B)  reduces (``at no cost'') to solving the two optimization problems:
\begin{equation}\label{eqArkadi616}
\begin{array}{ll}
(a)&\arg\min_{u\in Q_1}\left[\omega_x(u)+u^T[\zeta_x-\omega_x^\prime(\bar{x})]\right]\\
(b)&P_{\bar{Y}}(\zeta_Y):=\arg\min_{V\in Q_2}\left[\Ent(V)+\Tr\left(V\left[\zeta_Y-\Ent^\prime(\bar{Y})\right]\right)\right]\\
\end{array}
\end{equation}
with $\zeta_x\in \R^m=E_1$ and $\zeta_Y\in\mathcal{S}_n=E_2$ readily given by $\zeta$, specifically, $\zeta=(\Omega_x^{-2}\zeta_x,\Omega_Y^{-2}\zeta_Y)$. In the sequel, we assume that $(a)$ is easy to solve. The solution of $(b)$ can be written explicitly (see, e.g., \cite{Ben-tal_nemirovski_level}). Specifically, since $\bar{z}\in Q^o$, we have $\bar{Y}\in Q_2^o=\{Y\in\cS^n:Y\succ0,\Tr(Y)=1\}$, and the latter set clearly is the set of all matrices of the form
\begin{eqnarray}
 \label{definition_cH_of_V}
\cH(V)={\exp\{V\}\over\Tr(\exp\{V\})},\qquad V\in\cS^n.
\end{eqnarray}
\noindent Assuming that we have at our disposal a representation $\bar{Y}=\cH(\bar{V})$ with $\bar{V}\in\mathcal{S}_n$, the solution to $(b)$ is just $\cH(-\zeta_Y+\bar{V})$. In other words, when parameterizing points $Y\in Q_2^o$ according to $Y=\cH(V)$, prox-mapping (\ref{eqArkadi616}) becomes trivial - it reduces to a matrix addition. The $Y$-components of the points $w_t=(u_t,W_t)$ and $z_t=(x_t,Y_t)$ generated by the deterministic MP are of the form $W_t=P_{Y_{t-1}}(\Omega_Y^{2}[\eta_t]_Y)$ and $Y_t=P_{Y_{t-1}}(\Omega_Y^{2}[\zeta_t]_Y)$, where $\eta_t,\zeta_t\in \mathbb{R}^m\times \cS^n$ are given (see Algorithm \ref{alg:mirror_prox}).  When using parametric representations $W_t=\cH(U_t)$ and $Y_t=\cH(V_t)$, the matrices $U_t$ and $V_t$ are easy to update:  $U_t=V_{t-1}-\Omega_Y^{2}[\eta_t]_Y$ and $V_t=V_{t-1}-\Omega_Y^{2}[\zeta_t]_Y$, respectively, with $\eta_t$, $\zeta_t$ as defined in Algorithm \ref{alg:mirror_prox}. Thus, when representing $W_t$ and $Y_t$ by their ``matrix logarithms'' $U_t$ and $V_t$, it looks {\sl as if} the computational effort per step of MP as applied to (\ref{eqArkadi44}) were dominated by the necessity to resolve task (A), and in task (B) --- to solve the problem (\ref{eqArkadi616}.$a$) alone. This impression, however, is not fully true. Indeed, looking at (\ref{def:operatorF}), we observe that --- while computing the $Y$-component $F_Y(x)=-\cA(x)-B$ of $F$ at a point $z=(x,Y)$ needs the knowledge of $x$ only and is independent of how $Y$ is represented --- computing the $Y$-component $F_x(Y)=[\Tr(A_1Y);...;\Tr(A_mY)]+c$ of $F(z)$ seemingly requires the explicit representation of $Y$. This latter observation makes it necessary to solve explicitly problems (\ref{eqArkadi616}.$b$), or, which is the same, requires computation of the value of $\cH$ at a given point $V$. The related computational effort is $\mathcal{O}(n^3)$ (the arithmetic cost of an eigenvalue decomposition of $V$), which, depending on the problem's structure and sizes, can by far dominate all other computational expenses at an iteration. The goal of this paper is to demonstrate that one can avoid the explicit solution of ``troublemaking'' problems (\ref{eqArkadi616}.$b$), and use instead the easy-to-update ``logarithmic'' representations at the cost of a randomized computation of $F_x(\cdot)$. The idea of randomization is as follows: assume that we are given a ``matrix logarithm'' $V$ of $Y\in Q_2^o$, so that $Y=\cH(V)$. We need to compute a randomized estimate of the vector $F_x(Y)=\cA^\ast (Y)+c$, that is, of:
\[\cA^\ast (Y)=[\Tr(A_1Y);\ldots;\Tr(A_mY)]={1\over \Tr(\exp\{V\})}[\Tr(A_1\exp\{V\});\ldots;\Tr(A_m\exp\{V\})].\]
Imagine for a moment that we can multiply vectors by the matrix $\exp\{V/2\}$. Then, generating a sample $\xi$ of $N$ independent vectors $\xi^s\sim\cN(0,I_n)$, $1\leq s\leq N$,  and setting $\chi^s=\exp\{V/2\}\xi^s$, $s=1,2,...,N$, we have:
\begin{eqnarray*}
\bE\left\{{1\over N}\sum_{s=1}^N[[\chi^s]^TA_1\chi^s;\ldots;[\chi^s]^TA_m\chi^s]\right\}
&=&[\Tr(A_1\exp\{V\});...;\Tr(A_m\exp\{V\})],
\end{eqnarray*}
and:
\[\bE\left\{{1\over N}\sum_{s=1}^N[\chi^s]^T\chi^s\right\}=\Tr(\exp\{V\}),\]
so that we can use the random vector
\begin{equation}\label{estimate}
g_\xi(V):={\sum_{s=1}^N[[\chi^s]^TA_1\chi^s;...;[\chi^s]^TA_m\chi^s]\over \sum_{s=1}^N[\chi^s]^T\chi^s}+c,\qquad \chi^s:=\exp\{V/2\}\xi^s,
\end{equation}
as a random (biased!) estimate of $F_x(\cH(V))$. The last strategic question to be addressed is how indeed to compute, given $V$ and a vector $\xi$, the vector $\chi=\exp\{V/2\}\xi$. We propose to build a high accuracy approximation $\bar{\chi}$ to $\chi$ by setting
\begin{equation}\label{estimate1}
\bar{\chi}=\sum_{j=0}^J{1\over j!}(V/2)^j\xi
\end{equation}
 with $J$ large enough to guarantee a desired accuracy, and to compute the terms $\upsilon_j={1\over j!}(V/2)^j\xi$ by successive matrix-vector multiplications: $\upsilon_0=\xi$, $\upsilon_{j+1}={1\over 2(j+1)}V\upsilon_j$. We then merely replace in (\ref{estimate}) the vectors $\chi^s$ with their approximations $\bar{\chi}^s$, thus getting an estimate
\begin{equation}\label{estimate2}
\hat{g}_\xi(V):={\sum_{s=1}^N[[\bar{\chi}^s]^TA_1\bar{\chi}^s;...;[\bar{\chi}^s]^TA_m\bar{\chi}^s]\over \sum_{s=1}^N[\bar{\chi}^s]^T\bar{\chi}^s}+c
\end{equation}
of $g_\xi(V)$. We also set
\begin{equation}\label{hatcH}
\hat{\cH}_\xi(V)=\left[\sum_{s=1}^N[\bar{\chi}^s]^T\bar{\chi}^s\right]^{-1}\sum_{s=1}^N[\bar{\chi}^s[\bar{\chi}^s]^T];
\end{equation}
note that $\hat{\cH}_\xi(V)\in Q_2=\Delta_n^M$ can be considered as a random estimate of $\cH(V)$ (see (\ref{definition_cH_of_V})), and that $\hat{g}_\xi(V)=F_x(\hat{\cH}_\xi(V))$.

\subsubsection{Randomized algorithm}
Implementing the outlined randomization strategy with the setup presented in Section \ref{sect:setupMSP}, Algorithm \ref{alg:mirror_prox} becomes as follows:
\begin{myalgorithm}\label{MSPproblemMP}\label{alg:mirror_prox_adapted}{\rm [Randomized Mirror-Prox method applied to matrix saddle-point problem (\ref{eq:matrix_saddle-point})]\\
1: Choose the number of iterations $T$, the sample size $N$, and a sequence of positive integers $J_t$, $1\leq t\leq T$. Generate $2T$ independent samples $\xi_1,\xi_2,...,\xi_{2T}$, each of them comprised of $N$ independent realizations $\xi^s_t\sim\cN(0,I_n)$, $1\leq s\leq N$.\\
2: Set $x_0=x^{\omega_x}$ and let $V_0\in\mathcal{S}_n$ be the all zero matrix.\\
3: {\bf for} $1\leq t\leq T$ {\bf do}\\
4: Given $(x_{t-1},V_{t-1})$, choose (deterministic) $\gamma_t>0$ such that (cf. (\ref{eqArkadi6}), (\ref{eqArkadi66})):
\begin{equation} \label{gamma}
\gamma_t\leq \frac{1}{\sqrt{2}L}={1\over\sqrt{2}\Omega_x\Omega_Y\cL}.
\end{equation}
5: Compute the approximation
$$
\hat{F}_{\xi_{2t-1}}(x_{t-1},V_{t-1})=\left[\hat{g}_{\xi_{2t-1}}(V_{t-1}); -B-\cA(x_{t-1})\right]
$$
of $F(x_{t-1},\cH(V_{t-1}))=\left[\cA^\ast(\cH(V_{t-1}))+c;-B-\cA(x_{t-1})\right]$, where $\hat{g}_{\xi_{2t-1}}(\cdot)$ is as explained in (\ref{estimate2}), with $J_t$ in the role of $J$.
\\
6: Set
\begin{eqnarray}
\label{eq:update_wArkadi}
\bar{x}_t &=& \arg\min_{x\in Q_1} \left\{\left\langle \gamma_t \Omega_x^{2}\hat{g}_{\xi_{2t-1}}(V_{t-1}) -\omega_x^\prime(x_{t-1}),x\right\rangle+\omega_x(x)\right\}\cr
\bar{V}_t &=& V_{t-1}+\gamma_t\Omega_Y^2\left(B+\cA(x_{t-1})\right).
\end{eqnarray}
7: Compute the approximation $\hat{\cH}_t:=\hat{\cH}_{\xi_{2t}}(\bar{V}_t)$ of $\cH(\bar{V}_t)$ and the approximation
$$
\hat{F}_{\xi_{2t}}(\bar{x}_t,\bar{V}_t)=\left[\hat{g}_{\xi_{2t}}(\bar{V}_t);-B-\cA(\bar{x}_{t})\right]
$$
of $F(\bar{x}_t,\cH(\bar{V}_t))$, where $\hat{\cH}_{\xi_{2t}}(\bar{V}_t)$ and $\hat{g}_{\xi_{2t}}(\cdot)$ are as explained in (\ref{hatcH}) and  (\ref{estimate2}), respectively, and with $J_t$ in the role of $J$.\\
8: Set
\begin{eqnarray}
\label{eq:update_zArkadi}
x_t &=& \arg\min_{x\in Q_1} \left\{\langle \gamma_t \Omega_x^{2} \hat{g}_{\xi_{2t}}(\bar{V}_t)-\omega_x^\prime(x_{t-1}),x\rangle+\omega_x(x)\right\}\cr
V_t &=& V_{t-1}+\gamma_t\Omega_Y^2\left(B+\cA(\bar{x}_{t})\right).
\end{eqnarray}
9: {\bf end for}\\
10: Return $x^T:=\left(\sum_{t=1}^T\gamma_t\right)^{-1}\sum_{t=1}^T\gamma_t \bar{x}_t$.
}
\end{myalgorithm}

\subsubsection{Convergence and complexity analysis}
\paragraph{Regularity assumption and preliminaries.}
In order for Algorithm \ref{alg:mirror_prox_adapted} to be well behaved, we need certain additional assumption on $(E_1,\left\|\cdot\right\|_{x,*})$, specifically, the one of {\sl regularity} with certain parameter $\kappa=\kappa_{E_1}$. Instead of stating this notion here in full generality (this is done in Section \ref{sect:largedev} of the Appendix), let us just hint that the property has to do with ``good behavior'' of sums of martingale differences taking values in $E_1$ and list the regularity parameters for the most important, in regard to applications, pairs $(E_1,\left\|\cdot\right\|_{x,*})$. Specifically, denoting by $|\cdot|_p$ the spectral $\ell_p$-norm on the space $\cM^m$ of $m\times m$ matrices, that is, $|A|_p=\|\sigma(A)\|_p$, where $\sigma(A)$ is the vector of singular values of $A$, the following holds true (from now on, all $\mathcal{O}(1)$'s are appropriate absolute constants):
\begin{quote}
(!) {\sl If $1\leq p\leq2$ and either $(E_1,\left\|\cdot\right\|_x)=(\R^m,\left\|\cdot\right\|_p)$, or $(E_1,\left\|\cdot\right\|_x)=(\cM^m,|\cdot|_p)$, then the regularity parameter of $(E_1,\left\|\cdot\right\|_{x,*})$ is equal to 1 when $p=2$, is bounded from above by ${1\over p-1}$ when $p>1$, and is bounded from above by $\mathcal{O}(1)\ln(m+1)$.}
\end{quote}
From now on, if the opposite is not explicitly stated, it is assumed that $(E_1,\left\|\cdot\right\|_{x,*})$ is $\kappa$-regular.
An instrumental role in the convergence analysis of Algorithm \ref{MSPproblemMP} is played by the following fact (proved in Appendix).
\begin{prop}\label{prop:matrix_exp_approximation}
Let $(E_1,\left\|\cdot\right\|_{x,*})$ be $\kappa$-regular for some $\kappa$. With $F$ given by (\ref{def:operatorF}) and $g_\xi(\cdot)$ given by (\ref{estimate}), one has for every $V\in\cS^n$:
\begin{equation}\label{eqArkadi10}
\begin{array}{ll}
(a)&\|\bE_{\xi}\{g_{\xi}(V)\} -F_x(\cH(V))\|_{x,*}\leq \mathcal{O}(1)\cL\sqrt{\kappa}N^{-1}\\
(b)&\bE_{\xi}\left\{\exp\left\{{\sqrt{N}\|g_\xi(V)-\bE_{\xi}\{g_{\xi}(V)\}\|_{x,*}\over \mathcal{O}(1)\cL\sqrt{\kappa}}\right\}\right\}\leq\exp\{1\}.\\
\end{array}
\end{equation}
\end{prop}
Another component of our analysis is the following simple statement (recall that $\left\|\cdot\right\|_{Y,*}$ is the usual spectral norm on $\cS^n$):
\begin{prop}
\label{prop:matrix_exp}
Let $W\in \cS_n$ and $J\geq \exp(2)\|W\|_{ Y,*}$. Then,
\begin{eqnarray}
\label{eq:Taylor}
\left\|\exp\{W\}-\sum_{j=0}^J\frac{W^j}{j!}\right\|_{Y,*}\leq\exp\{-J\}.
\end{eqnarray}
\end{prop}
This result can be proved by applying the same arguments as in the proof of Lemma 6 in \cite{Arora_Kale_07}. 

\paragraph{Convergence analysis.} Proposition \ref{prop:matrix_exp} shows that in order to approximate the matrix exponent $\exp\{W\}$ by its Taylor polynomial within accuracy $\epsilon\ll1$, it suffices to take for the degree of the polynomial the number $J=\mathcal{O}(1)\ln(1/\epsilon)\|W\|_{Y,*}$, so that $J$ is ``nearly independent'' of $\epsilon$. Now, when $\epsilon$ is really small -- like $10^{-16}$ or even $10^{-100}$ -- any $\epsilon$-approximation of the matrix exponent is, ``for all practical purposes,'' the same as the matrix exponent itself. Assuming that the choice of $J$ indeed ensures ``really small'' inaccuracies in the approximation of the matrix exponent, we have all reasons to undertake a {\sl simplified} convergence analysis of Algorithm \ref{MSPproblemMP}, where we neglect the difference between the quantities $g_\xi(\cdot)$ as given by (\ref{estimate}) and their estimates $\hat{g}_\xi(\cdot)$ (defined in (\ref{estimate2})). Or, alternatively formulated: we analyze the idealized version of the algorithm with $g_\xi(\cdot)$ in place of $\hat{g}_\xi(\cdot)$.
\par
The convergence analysis of the idealized algorithm is as follows. Let $1\leq t\leq T$, and let $\gamma_t$ satisfy (\ref{gamma}). Note that in the notation of Algorithms \ref{generalMP} and \ref{alg:mirror_prox_adapted} and of definitions (\ref{martingale}) we have:
$$
\begin{array}{l}
z_{t-1}=(x_{t-1},\cH(V_{t-1})),\quad w_t=(\bar{x}_t,\cH(\bar{V}_t)),\\
\mu_{{z}_{t-1}}=\left[\bE_{\xi_{2t-1}}\left\{g_{\xi_{2t-1}}(V_{t-1})\big|\xi_{[2t-2]}\right\}-F_x(z_{t-1});0\right]=:[\mu^x_{z_{t-1}};0],\\
\mu_{{w}_{t}}=[\bE_{\xi_{2t}}\left\{g_{\xi_{2t}}(\bar{V}_{t})\big|\xi_{[2t-1]}\right\}-F_x(w_{t});0]=:[\mu^x_{w_{t}};0],\\
\sigma_{z_{t-1}}=\left[g_{\xi_{2t-1}}(V_{t-1})-\bE_{\xi_{2t-1}}\left\{g_{\xi_{2t-1}}(V_{t-1})\big|\xi_{[2t-2]}\right\};0\right]=:[\sigma^x_{z_{t-1}};0],\\
\sigma_{{w}_{t}}=\left[g_{\xi_{2t}}(\bar{V}_{t})-\bE\left\{g_{\xi_{2t}}(\bar{V}_{t})\big|\xi_{[2t-1]}\right\};0\right]=:[\sigma^x_{w_{t}};0].\\
\end{array}
$$
By (\ref{eqArkadi10}.$a$) combined with (\ref{conjnorm}), we obtain:
\begin{equation}\label{boundsmu}
\begin{array}{l}
\bE_{\xi_{2t-1}}\left\{\|\mu_{{z}_{t-1}}\|_*\big|\xi_{[2t-2]}\right\}=\Omega_x\bE_{\xi_{2t-1}}\left\{\|\mu^x_{{z}_{t-1}}\|_{x,*}\big|\xi_{[2t-2]}\right\}
\leq \mathcal{O}(1)\Omega_x\cL\sqrt{\kappa}N^{-1},\\
\bE_{\xi_{2t}}\left\{\|\mu_{{w}_{t}}\|_*\big|\xi_{[2t-1]}\right\}=\Omega_x\bE_{\xi_{2t}}\left\{\|\mu^x_{{w}_{t}}\|_{x,*}\big|\xi_{[2t-1]}\right\}
\leq \mathcal{O}(1)\Omega_x\cL\sqrt{\kappa}N^{-1},\\
\end{array}
\end{equation}
whence also:
\begin{equation}\label{also1}
\begin{array}{l}
\bE_{\xi_{[2T]}}\left\{\sum_{t=1}^T\gamma_t\|\mu_{{w}_t}\|_*\right\}=\Omega_x\bE_{\xi_{[2T]}}\left\{\sum_{t=1}^T\gamma_t\|\mu^x_{{w}_t}\|_{x,*}\right\}
\leq \mathcal{O}(1)\Omega_x\cL\sqrt{\kappa}N^{-1}\sum_{t=1}^T\gamma_t,\\
\begin{array}{ll}
\bE_{\xi_{[2T]}}\left\{\sum_{t=1}^T\gamma_t^2\|\mu_{{w}_t}-\mu_{z_{t-1}}\|_*^2\right\}&=\Omega_x^2
\bE_{\xi_{[2T]}}\left\{\sum_{t=1}^T\gamma_t^2\|\mu_{{w}_t}-\mu_{z_{t-1}}\|_{x,*}^2\right\}\\
&\leq \mathcal{O}(1)\Omega_x^2\cL^2\kappa
N^{-2}\sum_{t=1}^T\gamma_t^2.\\
\end{array}
\end{array}
\end{equation}
Further, inequality (\ref{eqArkadi10}.$b$) implies:
\begin{equation}\label{boundssigma}
\bE_{\xi_{2t-1}}\left\{\|\sigma^x_{z_{t-1}}\|_{x,*}^2\big|\xi_{[2t-2]}\right\}\leq
\mathcal{O}(1)\cL^2\kappa/N,\qquad
\bE_{\xi_{2t}}\left\{\|\sigma^x_{w_{t}}\|_{x,*}^2\big|\xi_{[2t-1]}\right\}
\leq
\mathcal{O}(1)\cL^2\kappa/N.
\end{equation}
Moreover, by the definition of $\sigma^x_{z_{t-1}}$ and $\sigma^x_{w_t}$, we have:
\begin{equation}\label{clearly}
\bE_{\xi_{2t-1}}\{\sigma^x_{z_{t-1}}\big|\xi_{[2t-2]}\}=0,\qquad \bE_{\xi_{2t}}\{\sigma^x_{w_{t}}\big|\xi_{[2t-1]}\}=0.
\end{equation}
Since $(E_1,\left\|\cdot\right\|_{x,*})$ is $\kappa$-regular, relations (\ref{boundssigma}) and (\ref{clearly}) imply by
Proposition 3 of \cite{Nemirovski_Regular_Banach_Spaces} that:
$$
\bE_{\xi_{[2T]}}\left\{\left\|\sum_{t=1}^T\gamma_t\sigma^x_{w_t}\right\|_{x,*}^2\right\}\leq \mathcal{O}(1)\kappa^2\mathcal{L}^2 N^{-1} \sum_{t=1}^T\gamma_t^2,$$
which results in:
\begin{equation}\label{implythat}
\begin{array}{rcl}
\bE_{\xi_{[2T]}}\left\{\left\|\sum_{t=1}^T\gamma_t\sigma_{w_t}\right\|_*\right\}&=&\Omega_x\bE_{\xi_{[2T]}}\left\{\left\|\sum_{t=1}^T\gamma_t\sigma_{w_t}^x\right\|_{x,*}\right\}
\leq\Omega_x\sqrt{\bE_{\xi_{[2T]}}\left\{\left\|\sum_{t=1}^T\gamma_t\sigma^x_{w_t}\right\|_{x,*}^2\right\}}\\
&\leq& \mathcal{O}(1)\Omega_x\cL\kappa N^{-1/2}\sqrt{\sum_{t=1}^T\gamma_t^2}.\\
\end{array}
\end{equation}
Besides this, (\ref{boundssigma}) implies that:
\begin{eqnarray}\label{implythatagain}
\bE_{\xi_{[2T]}}\left\{\sum_{t=1}^T\gamma_t^2\|\sigma_{{w}_t}-\sigma_{{z}_{t-1}}\|_*^2\right\}&=&\Omega_x^2
\bE_{\xi_{[2T]}}\left\{\sum_{t=1}^T\gamma_t^2\|\sigma^x_{{w}_t}-\sigma^x_{{z}_{t-1}}\|_{x,*}^2\right\}\cr
&\leq& \mathcal{O}(1)\Omega_x^2\cL^2\kappa N^{-1}\sum_{t=1}^T\gamma_t^2.
\end{eqnarray}
Combining (\ref{also1}), (\ref{implythat}), (\ref{implythatagain}), and taking into account that with our setup $\Omega=\sqrt{2}$ and $D\leq 2\Omega$, we conclude from Theorem \ref{thm:convergence} that:
\begin{equation*}
\begin{array}{rcl}
\bE\left\{\overline{\phi}(x^T)-\min_{x\in Q_1}\overline{\phi}(x)\right\}
&\leq&\mathcal{O}(1){1+\Omega_x\cL\sqrt{\kappa}\left[N^{-1}\sum_{t=1}^T\gamma_t+\sqrt{\kappa}N^{-1/2}\sqrt{\sum_{t=1}^T\gamma_t^2}
+\Omega_x\cL\sqrt{\kappa}N^{-1}\sum_{t=1}^T\gamma_t^2\right]\over \sum_{t=1}^T\gamma_t}
\end{array},
\end{equation*}
where $\phi(x,Y)$ is the cost function of the saddle point problem (\ref{eqArkadi44}), $\overline{\phi}(x)=\max\limits_{Y\in\Delta_n^M}\phi(x,Y)$ is the objective in the problem of interest (\ref{eqArkadi4}), and $\Opt=\min\limits_{x\in Q_1}\overline{\phi}(x)$ is the saddle point value in (\ref{eqArkadi44}), or, equivalently, the optimal value in (\ref{eqArkadi4}).
\par
Optimizing the resulting efficiency estimate in the stepsizes $\gamma_t$ satisfying (\ref{gamma}), it is immediately seen that with
\begin{equation}\label{parametersfinal}
N=\hbox{\rm floor}\left({T\kappa\over 2\Omega_Y^2}\right)=\hbox{\rm floor}\left({T\kappa\over 4\ln(n)}\right),\qquad \gamma_t={
1\over\sqrt{2}\Omega_x\Omega_Y\cL}={1\over 2\sqrt{\ln(n)}\Omega_x\cL},\ 1\leq t\leq T,
\end{equation}
the above inequality implies:
\begin{equation}\label{finalbound}
\bE\left\{\overline{\phi}(x^T)-\min_{x\in Q_1}\overline{\phi}(x)\right\}
\leq \mathcal{O}(1){\Omega_x\cL\over T}\left[{\ln(n)\over\sqrt{\kappa}}+\sqrt{\kappa\ln(n)}\right].
\end{equation}

\section{An application: minimizing the maximal eigenvalue of a convex combination of symmetric matrices}
\label{sec:applications}

Consider the special case of problem (\ref{eqArkadi4}) where $Q_1$ is the standard simplex in $\R^m$:
$$
Q_1=\Delta_m:=\left\{x\in\R^m:x\geq0,\sum_{j=1}^mx_j=1\right\}\eqno{[m\geq2]}
$$
and $B=0$, $c=0$, so that the problem is
\begin{equation}\label{problem:min_max_eigenvalue}
\Opt=\min\limits_{x\in\Delta_m}\lambda_{\max}(\cA(x)),\qquad \cA(x)=\sum_{j=1}^mx_jA_j.
\end{equation}
In other words, we want to minimize the largest eigenvalue of a convex combination $\sum_jx_jA_j$ of given symmetric $m\times m$ matrices. Note that the problem of minimizing the maximal eigenvalue of $B+\cA(x)$ over $Q_1$ reduces to (\ref{problem:min_max_eigenvalue}) by replacing the matrices $A_j$ with $A_j+B$. The operator (\ref{def:operatorF}) associated with the problem is
\begin{equation}\label{assocF}
F(x,Y)=\left[F_x(Y);F_y(x)\right]:=\left[[\Tr(A_1Y);...;\Tr(A_mY)];-\cA(x)\right].
\end{equation}
\par
We equip $E_1=\R^m$ with the norm $\|\cdot\|_x=\|\cdot\|_1$; the conjugate norm is $\|\cdot\|_{x,*}=\|\cdot\|_\infty$. It is well known that $(\R^m,\|\cdot\|_\infty)$ is $\kappa$-regular with $\kappa=\mathcal{O}(1)\ln(m)$ (see, e.g., Example 2.1 in \cite{Nemirovski_Regular_Banach_Spaces}). Note that this choice of $\|\cdot\|_1$ results in
\begin{equation}\label{cLis}
\cL=\max\limits_{1\leq j\leq m} \|A_j\|_{Y,*}
\end{equation}
(see (\ref{eqArkadi6})), where $\|A\|_{Y,*}$ is the spectral norm of a matrix $A$.
\par
We equip $Q_1$ with the d.-g.f. function:
$$
\omega_x(x)=\Ent(x):=\sum_{j=1}^mx_j\ln(x_j),
$$
which, as it is well known (and immediately seen), is compatible with $\|\cdot\|_1$. The associated entities are
\begin{equation}\label{Omegaxetc}
Q_1^o=\left\{x\in\R^m:x>0,\sum_{j=1}^mx_j=1\right\},\qquad x^\omega=\left[{1\over m};...;{1\over m}\right],\qquad\Omega_x=\sqrt{2\ln(m)},
\end{equation}
and the prox-mapping is given by an explicit formula:
$$
\left(\arg\min_{w\in Q_1}\left\{\omega_x(w)+\langle \xi-\omega_x^\prime(x),w\rangle\right\}\right)_j={\exp\{\ln(x_j)-\xi_j\}\over \sum_{\ell=1}^m\exp\{\ln(x_\ell)-\xi_\ell\}},\qquad 1\leq j\leq m,
$$
see Section \ref{subs:strategy}.
\par
Let us solve (\ref{problem:min_max_eigenvalue}) by $T$-step Algorithm \ref{MSPproblemMP} associated with the outlined setup and the parameters $N$ and $\gamma$ chosen according to (\ref{parametersfinal}), that is, as:
\begin{equation}\label{parametersagain}
N=\mathcal{O}(1){\ln(m)\over\ln(n)}T, \qquad \gamma_t\equiv\gamma={1\over 2\cL \sqrt{2\ln(n)\ln(m)}}.
\end{equation}
The efficiency estimate (\ref{finalbound}) now reads:
\begin{equation}\label{nowreads}
\bE\left\{\lambda_{\max}(\cA(x^T))-\Opt\right\}\leq \mathcal{O}(1)\left(\ln(n)+\ln(m)\sqrt{\ln(n)}\right)\cL/T.
\end{equation}
\paragraph{Choosing truncation levels $J_t$.} Let us specify the ``truncation levels'' $J_t$ for $1\leq t\leq T$. In view of (\ref{eq:update_wArkadi}), (\ref{eq:update_zArkadi}), (\ref{cLis}) and taking into account that $V_0=0$, $\Omega_Y=\sqrt{2\ln(n)}$, and $\Omega_x=\sqrt{2\ln(m)}$, we conclude that:
$$
\|V_t\|_{Y,*}\leq \mathcal{O}(1)\sqrt{\ln(n)/\ln(m)}t,\qquad \|\bar{V}_t\|_{Y,*}\leq \mathcal{O}(1)\sqrt{\ln(n)/\ln(m)}t.
$$
From Proposition \ref{prop:matrix_exp}, we deduce that the matrix exponentials we need to use can be approximated with accuracy $\delta\ll 1$ by a truncated Taylor series with $J_t=\mathcal{O}(1)\sqrt{\ln(n)/\ln(m)}\ln(1/\delta)t$ terms. Specifying $\delta$ as, say, machine accuracy, we see that ``for all practical purposes'' it suffices to take
\begin{equation}\label{eqJt}
J_t=\mathcal{O}(1) \sqrt{\ln(n)/\ln(m)}t,\qquad t\geq1,
\end{equation}
with a moderate absolute constant $\mathcal{O}(1)$.
\paragraph{Overall complexity.} Assume that we want to solve (\ref{problem:min_max_eigenvalue}) within accuracy $\epsilon$ in terms of the objective. This task is typically unreachable with a randomized algorithm. Instead, we need to content ourselves with a procedure returning an $\epsilon$-solution with a prescribed probability of at least $1-\beta$, where $0<\beta\ll 1$. To build such a procedure, we can specify $T=T(\epsilon)$ in such a way that the right hand side in (\ref{nowreads}) is at most $\epsilon/4$. We run the above $T(\epsilon)$-step algorithm $k$ times,  each time computing an accurate, within the margin $\epsilon/2$, estimate of the value $\lambda_{\max}(\cA(x^{T,i}))$ of the objective at the corresponding output $x^{T,i}$, $1\leq i\leq k$, and then select among the $k$ outputs $x^{T,1},...,x^{T,k}$ the one with the smallest estimate of the objective value. Since with our choice of $T(\epsilon)$ we have $\Prob\{\lambda_{\max}(x^{T,i})-\Opt>\epsilon/2\}\leq {1\over 2}$ and $x^{T,1}$,...,$x^{T,k}$ are independent, this procedure yields an $\epsilon$-solution to the problem of interest with a probability of at least $1-\beta$ for a ``small'' $k=\mathcal{O}(1)\ln(1/\beta)$.
\par
Now, let us evaluate the computational complexity of a single $T(\epsilon)$-step run of Algorithm \ref{MSPproblemMP}. Assume that every matrix $A_i$ has at most $S$ nonzero entries. We assume that $mS\geq n^2$, meaning that the matrices $\cA(x)$ can be fully dense. In order to avoid intricate expressions, we omit in the sequel all factors that are logarithmic in $m$, $n$ and $1/\beta$ (in particular, all absolute constant factors) and write down the statement ``$P$ is, within logarithmic factors,  bounded from above by $Q$'' as $P\lesssim Q$. We also write $P\sim Q$ when both $P\lesssim Q$ and $Q\lesssim P$.  Finally, we set $\nu=\epsilon/\cL$; this quantity can be naturally interpreted as the relative accuracy of an $\epsilon$-solution. To establish the complexity of our procedure, note the following.
\begin{itemize}
\item[(A) ] By (\ref{nowreads}), the required number of steps $T=T(\epsilon)$ admits the bound $T\lesssim1/\nu$, whence, by (\ref{parametersagain}), $N\lesssim 1/\nu$.
\item[(B) ] As it is immediately seen, when $mS\geq n^2$, the computational effort at step $t\leq T$ of the algorithm is, within factor $\mathcal{O}(1)$,  dominated by the necessity
\begin{enumerate}
\item to compute $\cA(x)$ at a given point $x$, ($\lesssim mS$ arithmetic operations (a.o.));
\item to generate $N$ samples $\xi_t^s\sim\cN(0,I_n)$ with $1\leq s\leq N$, (totally $\lesssim nN$ a.o.);
\item to compute for every $s\leq N$ the vectors $\bar{\chi}_t^s=\sum_{j=0}^{J_t}(V/2)^j\xi_t^s/j!$, where $V\in\cS^n$ is a given matrix (totally $\lesssim Ntn^2$ a.o., see (\ref{eqJt}));
\item to build the matrix $H=\left[\sum_{s=1}^N[\bar{\chi}_t^s]^T\bar{\chi}_t^s\right]^{-1}\sum_{s=1}^N\bar{\chi}_t^s[\bar{\chi}_t^s]^T$ ( $\lesssim Nn^2$ a.o.);
\item  to compute the vector $[\Tr(HA_1);\Tr(HA_2);...;\Tr(HA_n)]$ ($\lesssim mS$ a.o.).
\end{enumerate}
\end{itemize}
We see that the complexity of step $t$ is $\lesssim Ntn^2+mS$ a.o., implying that the overall complexity of a single run of the algorithm is $\lesssim NT^2n^2+TmS\lesssim n^2/\nu^3+mS/\nu$ a.o. We then should compute the value of the objective at the resulting approximate solution, that is, the maximal eigenvalue of a symmetric matrix with the spectral norm not exceeding $\cL$. For our purposes, it suffices to approximate this value $(1-\beta/k)$-reliably within accuracy $\mathcal{O}(\epsilon)$, which can be done by the Power method at the cost of $\lesssim n^2/\nu$ a.o. Finally, we should repeat this procedure $\mathcal{O}(1)\ln(1/\beta)$ times. Omitting constants and factors logarithmic in $m,$ $n,$ and $1/\beta$, our randomized algorithm yields an $(1-\beta)$-reliable $\epsilon$-solution to {\rm (\ref{problem:min_max_eigenvalue})} at the cost of 
$$
\cC_{\hbox{\scriptsize SMP}}={n^2\over\nu^3}+{mS\over\nu}\text{ a.o.}\eqno{[\nu=\epsilon/\cL].}
$$
\paragraph{Discussion.} Let us compare the complexity of our algorithm with those of its existing competitors. To the best of our knowledge, the best existing complexity bounds for large-scale problems (\ref{problem:min_max_eigenvalue}) are as follows (we again skip logarithmic factors):
\begin{itemize}
\item The complexity for {\sl Interior Point methods} without any assumptions on $A_j$ aside of their sparsity is
\[\cC_{\hbox{\scriptsize IPM}}=\sqrt{\max[n,m]}[n^3+m^3+m^2n^2]\ln(1/\nu)\text{ a.o.}\]
\item {\sl Advanced deterministic first-order algorithms}, like Nesterov's Smoothing \cite{nesterov:coreDP73/2004} or deterministic Mirror-Prox, require
\[\cC_{\hbox{\scriptsize FOM}}={n^3+m S\over \nu}\text{ a.o.}\]
\item  We can also consider minimizing the original objective function $x\mapsto \lambda_{\max}(\sum_jx_jA_j)$ over the standard simplex using a {\sl``slightly randomized''  Mirror Descent method} \cite{dAspremont_Subsampling}. This method requires
\[\cC_{\hbox{\scriptsize MD}}={n^2\over\nu^{5/2}}+{mS\over\nu^2}\text{ a.o.}\]
The iteration count in this method is $\sim 1/\nu^2$. The computational effort per iteration reduces to assembling $\cA(x)$ at a given point ($\sim mS$ a.o.) and computing an $\epsilon$-subgradient of the objective and an $\epsilon$-approximation of the value of the objective at $x$ by applying the Power method to the matrix $\cA(x)$ in order to approximate its maximal singular value and leading eigenvector. With a straightforward implementation of the Power method this task requires $\sim n^2/\nu$ a.o., and with an advanced implementation $\sim n^2/\sqrt{\nu}$ a.o. only.
\end{itemize}
It turns out that {\sl there exists a meaningful range of values of $m$, $n$, $S$, and $\nu$ where our stochastic algorithm significantly outperforms the outlined competitors.} For example, consider the case when $n$ is large, and assume that we have for some $0<\kappa<1/4$:
$$
mS\sim n^\beta \text{\ with\ } 2+2\kappa\leq \beta\leq 3-2\kappa,\,\,n^{\max[2-\beta,-1/2]+\kappa}\leq\nu\leq n^{1-\beta/2}
$$ (note that the outlined range of values of $\nu$ is nonempty; e.g., this range is $n^{-1/2+\kappa}\leq \nu\leq n^{-1/4}$ with $\beta=2.5$ ). It is immediately seen that in the case in question we have $\cC_{\hbox{\scriptsize SMP}}\sim n^2/\nu^3$ and:
$$
{\cC_{\hbox{\scriptsize IPM}}\over \cC_{\hbox{\scriptsize SMP}}}\gtrsim n^{3\kappa},\,\, {\cC_{\hbox{\scriptsize FOM}}\over\cC_{\hbox{\scriptsize SMP}}}
\gtrsim n^{2\kappa},\,\,{\cC_{\hbox{\scriptsize MD}}\over \cC_{\hbox{\scriptsize SMP}}}\geq n^\kappa,
$$
that is, our algorithm progressively outperforms its competitors as the sizes grow.

\section{Numerical experiments}

We consider randomly generated instances of the problem
\begin{equation}\label{problem:min_max_eigenvalue_random}
\Opt=\min\limits_{x\in\Delta_m}\lambda_{\max}(\cA(x)),\qquad \cA(x)=\sum_{j=1}^mx_jA_j, \qquad A_j=j^{3/2}C_j,
\end{equation}
where $C_j$ is a sparse symmetric random $n\times n$-matrix and $j=1,\ldots,m$, i.e., we are confronted with instances of problem (\ref{problem:min_max_eigenvalue}) that we studied in the last section. We solve these problem instances up to a (relative) accuracy of $\delta:=\epsilon\mathcal{L}$, where $0<\epsilon<1$ is the target accuracy and $\mathcal{L}$ is defined as in (\ref{cLis}). In all the numerical experiments that we perform, the target accuracy $\epsilon$ is set to $2\cdot 10^{-3}$.

\begin{table}
\label{table:N}
\begin{center}
\begin{tabular}[]{|c|ccc|ccc|c|}
\hline
 & \multicolumn{3}{c|}{CPU time [sec]} &\multicolumn{3}{c|}{number of iterations} &  CPU time per iteration \cr
 $N$ & mean & std & $95$\% conf & mean & std & $95$\% conf & [sec/iteration] \cr
\hline
$1$ & $66$ & $9$ & $[60,72]$  & $2948$ & $327$ & $[2746,3151]$ & $0.0224$ \cr

$3$ & $90$ & $14$ & $[81,98]$  & $2970$ & $343$ & $[2757,3183]$ & $0.0302$ \cr

$5$ & $86$ & $11$ & $[79,93]$ & $2900$ & $271$ & $[2732,3068]$ & $0.0298$ \cr

$10$ & $87$ & $12$ & $[80,94]$& $2860$ &$207$ & $[2732,2988]$ & $0.0305$ \cr

$50$ & $92$ & $7$ & $[87,96]$ & $2840$ &$232$ & $[2696,2984]$ & $0.0323$ \cr

$100$ & $98$ & $9$ & $[93,104]$& $2850$ &$207$ & $[2722,2978]$& $0.0344$  \cr

$500$ & $141$ & $15$ & $[131,150]$& $2860$ &$222$ & $[2722,2998]$ & $0.0491$ \cr

$1000$ & $178$ & $18$ & $[167,189]$& $2860$ &$222$ & $[2722,2998]$& $0.0622$  \cr

$5000$ & $533$ & $47$ & $[504,562]$ & $2877$ &$204$ & $[2750,3003]$ & $0.1851$ \cr
\hline
\end{tabular}\\
\end{center}
\caption{CPU time (mean, standard deviation, and $95$\% confidence interval), number of iterations (mean, standard deviation, and $95$\% confidence interval), and average CPU time per iteration required by the stochastic Mirror-Prox method for solving random instances of problem (\ref{problem:min_max_eigenvalue_random}) with parameter values $n=100$, $m=100$, $\epsilon=0.002$, and for different samples sizes $N$. The matrices $A_j$ have a joint sparsity pattern and, in expectation, $10$\% of the entries are non-zero.} 
\end{table}

 We implement Algorithm \ref{alg:mirror_prox_adapted} with constant step-sizes $\gamma_t=\gamma$, $t \geq 1$, and $\gamma$ has the form described in (\ref{parametersagain}). Given a matrix $W\in\mathcal{S}_n$, we choose the truncation level $J_W$ of the matrix exponential Taylor series approximation according to the following formula (compare with Proposition \ref{prop:matrix_exp}):
\[J=\left\lfloor \max\left\{ \log(1/\rho),\exp(1)\left\|W\right\|_{(\infty)}\right\}\right\rfloor,\qquad \rho:= 10^{-3}.\]
Note that this setting slightly deviates from the truncation level derived in Proposition \ref{prop:matrix_exp}. The $\infty$-norm of $W$ is computed approximately using the Power method. In accordance to (\ref{parametersagain}) and (\ref{nowreads}), we need to choose the sample size $N$ as:
\begin{eqnarray}
 \label{eq:N_numerical_results}
N=\frac{\mathcal{O}(1)\ln^2(m)}{\epsilon\sqrt{\ln(n)}}.
\end{eqnarray}
In Table 1, we give the CPU time (mean, standard deviation, and corresponding $95$\% confidence interval), the number of iterations needed to find a solution with relative accuracy $\epsilon\mathcal{L}$ (mean, standard deviation, and corresponding $95$\% confidence interval), and the average CPU time per iteration for different samples sizes $N$. The matrices $A_j$ are all of size $100\times100$, they follow the same sparsity pattern, and, on average, $10$\% of the entries are different from zero. In total, we have a hundred matrices $A_j$. All the numerical results that we present in this paper are averaged over ten runs and are obtained on a computer with 24 processors, each of them with $2.67$ GHz, and with $96$ GB of RAM. We observe that the smaller the sample size the lower the CPU time that is required to approximately solve the problem instances. Surprisingly, we can choose a very small sample size without sacrificing too many iterations. Let us illustrate this observation with an example. According to (\ref{eq:N_numerical_results}) and with an absolute constant of $1$, we are supposed to choose $N$ as about $5000$. With this parameter choice, we need an average CPU time of $533$ seconds. Using only one sample for each matrix exponential approximation, we observe that we can reduce the average CPU time by $87.6\%$ (with a slight increase of $2.5\%$ in the average number of iterations). For the subsequent tests, we will thus choose a sample size that deviates from its theoretical value and use only one sample for every matrix exponential approximation. 

\begin{sidewaystable}

\begin{center}
\begin{tabular}[]{|c|cc|cc|cc|ccc|}
\multicolumn{10}{c}{\textbf{CPU time (mean [sec], standard deviation [sec])}}\cr
\hline
& \multicolumn{2}{c|}{Mirror-Descent} &\multicolumn{2}{c|}{det Mirror-Prox} &\multicolumn{2}{c|}{stoch Mirror-Prox}&\multicolumn{3}{c|}{ratios CPU time}\cr
$(n,S)$ & mean & std &  mean & std & mean & std  & $\frac{\text{det MP}}{\text{MD}}$ & $\frac{\text{stoch MP}}{\text{MD}}$& $\frac{\text{stoch MP}}{\text{det MP}}$ \cr
\hline
$(100,955)$ & $307$ & $47$ & $128$ & $16$  & $71$ & $6$   & $0.42$& $0.23$ & $0.55$\cr
$(200,3813)$ & $766$ & $79$ & $307$ &$15$    & $237$ & $17$ & $0.40$&  $0.31$& $0.77$ \cr
$(400,15255)$ & $2522$ & $120$ & $1101$ &$42$  & $744$ & $39$ & $0.44$&  $0.30$& $0.68$ \cr
$(800,60971)$ & $10262$ & $746$ & $4983$ &$126$  & $2814$ & $74$&  $0.49$&  $0.27$& $0.56$\cr
\hline
\end{tabular}\\
\vspace{0.5cm}
\begin{tabular}[]{|c|ccc|ccc|ccc|ccc|}
\multicolumn{13}{c}{\textbf{number of iterations (mean, standard deviation) and average CPU time per iteration [sec / iteration]}}\cr
\hline
& \multicolumn{3}{c|}{Mirror-Descent} &\multicolumn{3}{c|}{det Mirror-Prox} &\multicolumn{3}{c|}{stoch Mirror-Prox}&\multicolumn{3}{c|}{ratios CPU time / iteration}\cr
$(n,S)$ & mean & std & $\frac{\text{CPU time}}{\text{iteration}}$ & mean & std & 
$\frac{\text{CPU time}}{\text{iteration}}$ & mean & std & $\frac{\text{CPU time}}{\text{iteration}}$  & $\frac{\text{det MP}}{\text{MD}}$ & $\frac{\text{stoch MP}}{\text{MD}}$& $\frac{\text{stoch MP}}{\text{det MP}}$ \cr
\hline
$(100,955)$ & $21504$ & $3287$ & $0.0143$ & $3120$ & $349$ & $0.0411$  & $3000$ & $313$ & $0.0235$ & $2.88$ & $1.65$ & $0.57$ \cr

$(200,3813)$ & $21388$ & $1858$ & $0.0358$& $2700$ &$141$ & $0.1137$   & $2740$ & $171$ & $0.0866$ & $3.17$& $2.42$ & $0.76$ \cr

$(400,15255)$ & $22293$ & $924$ & $0.1131$ & $2680$ &$103$ & $0.4109$  & $2620$ & $103$ & $0.2841$ & $3.63$ &$2.51$& $0.69$\cr

$(800,60971)$ & $22327$ & $445$ & $0.4596$ & $2740$ &$52$ & $1.8187$  & $2600$ & $47$ & $1.0822$ & $3.96$& $2.35$ & $0.60$ \cr
\hline
\end{tabular}\\
\vspace{0.5cm}
\begin{tabular}[]{c}
\textbf{Truncation level $J$ of Taylor series approximation (only stochastic Mirror-Prox)}\cr
\end{tabular}\\
\begin{tabular}[]{|c|cc|}
\hline
$(n,S)$ & mean & std \cr
\hline
$(100,955)$ & $9$ & $<1$  \cr

$(200,3813)$ & $9$ & $<0.5$ \cr

$(400,15255)$ & $10$ & $<0.5$  \cr

$(800,60971)$ & $10$ & $<0.5$   \cr
\hline
\end{tabular}\\
\end{center}
\label{table:comparison}
\caption{CPU time (mean and standard deviation), number of iterations (mean and standard deviation), and average CPU time per iteration needed by the Mirror-Descent (MD), the deterministic Mirror-Prox (det MP), and the stochastic Mirror-Prox (stoch MP) method for solving random instances of problem (\ref{problem:min_max_eigenvalue_random}) with parameter values $m=100$ and $\epsilon=0.002$, with $S$ non-zero entries, and for different values of the matrix size $n$. The performance ratios express the CPU time (CPU time per iteration) required by "method A" as percentage of the corresponding quantity used by "method B". The stochastic Mirror-Prox method is implemented with $N=1$, and the used truncation levels $J$ are shown in the table at the bottom.} 
\end{sidewaystable}

Given a pair $(\bar{x},\bar{Y})$ of a primal and a dual feasible solution, we can compute the corresponding duality gap
\begin{eqnarray}
\label{eq:stopping_criterion}
\lambda_{\max}(\cA(\bar{x}))- \min_{x\in\Delta_m}\left\langle \cA(x), \bar{Y} \right\rangle_F,
\end{eqnarray}
which we use as stopping criterion for our algorithm and which we check at every $100$th iteration of the method. The first term is approximated by an adapted version of the Power method and the second term is simply $\min\{\langle A_j,\bar{Y}\rangle_F:1\leq j\leq m\}$. As the Power method typically returns a lower bound on the eigenvalue of largest absolute value, we recompute the duality gap using the Matlab built-in functions \textit{max()} and \textit{eig()} when the first approximation obtained by our version of the Power method yields to a value that is smaller than $\epsilon\mathcal{L}$. We denote by $\hat{\mathcal{H}}(V)$ the approximation of $\exp(V)/\Tr(\exp(V))$ by the truncated Taylor development. The pair $(\bar{x},\bar{Y})$ considered at iteration $t$ is the average
\[
\frac{1}{\sum_{\tau=1}^t\gamma_\tau}\sum_{\tau=1}^t\gamma_\tau(\bar{x}_\tau,\hat{\mathcal{H}}(\bar
V_\tau))=\frac{1}{t}\sum_{\tau=1}^t(\bar{x}_\tau,\hat{\mathcal{H}}(\bar
V_\tau)),
\]
where $\bar{x}_\tau$ and $\bar V_\tau$ are defined in Algorithm
\ref{alg:mirror_prox_adapted}, equations (\ref{eq:update_wArkadi}). In principle, the criterion (\ref{eq:stopping_criterion}) gives theoretically a desirable solution only if we use exact scaled exponentials instead of $\hat{\mathcal{H}}(V_\tau)$. Nevertheless, $\hat{\mathcal{H}}(V)$ is in the matrix simplex by construction, and the number of terms we use in the Taylor exponential is large enough to justify a very accurate approximation, so
that $\hat Y$ can be considered as an adequate approximate solution to our problem. 

In Table 2, we compare the performance of our randomized version of the Mirror-Prox method with the efficiency of its deterministic counterpart and of the Mirror-Descent scheme for random problem instances (\ref{problem:min_max_eigenvalue_random}). As before, we have a hundred matrices $A_j$, but this time their size is varying. They are sparse with a joint sparsity pattern and with $S$ non-zero values; the values for $S$ can be found in Table \ref{table:comparison}. In this table, we show the CPU time (mean and standard deviation), the number of iterations required to find a solution with accuracy $\epsilon\mathcal{L}$ (mean and standard deviation), and the average CPU time per iteration. Moreover, we express the average (total) CPU time and the average CPU time per iteration of the Mirror-Descent method (deterministic Mirror-Prox) in percentage of the stochastic and the deterministic Mirror-Prox method (stochastic Mirror-Prox). We observe that the stochastic Mirror-Prox method has an average CPU time that corresponds to $23$ to $31\%$ of the running of the Mirror-Descent scheme and to $55$ to $77$\% of the CPU time required by the deterministic Mirror-Prox method for problem instances involving matrices of size $100\times 100$ up to size $800\times 800$. For the stochastic Mirror-Prox method, we also give the truncation levels $J$, which are nine or ten on average for these problem instances.  

\textbf{Acknowledgments:} We would like to thank Hans-Jakob L\"uthi for having initiated this collarobation and for many helpful discussions. This research was partially funded NSF Grant DMS-0914785 and by Swiss National Science Foundation (SNF), project no 200021-129743, ``First-order methods for Semidefinite Optimization''.

\bibliographystyle{amsalpha}
\bibliography{Jordan}

\appendix

\section{Large deviations of random sums}\label{sect:largedev}

Let $E$ be a Euclidean space with inner product $\langle \cdot, \cdot
\rangle$ and let $\kappa\geq 1$. We endow $E$ with a norm denoted by
$\lVert\cdot\rVert$, which might differ from the one associated with
this inner product.

We say that the space $(E,\lVert\cdot\rVert)$ is $\kappa$-smooth, if
the function $p(\cdot):=\lVert \cdot \rVert^2$ is continuously
differentiable on $E$ and
\[p(x+y)\leq p(x)+\langle p'(x),y\rangle +\kappa p(y) \qquad \forall\ x,y\in E.\]

Both the space $(E,\lVert\cdot\rVert)$ and the norm
$\lVert\cdot\rVert$ are called $\kappa$-regular, if there exist
$\kappa_+\in [1,\kappa]$ and a norm $\lVert\cdot\rVert_+$ on $E$ with
the following two characteristics:
\begin{itemize}
\item The space $(E,\lVert\cdot\rVert_+)$ is $\kappa_+$-smooth. \item
The norm $\lVert\cdot\rVert_+$ is $\kappa/\kappa_+$-compatible with
$\lVert\cdot\rVert$, that is:
\[\lVert x \rVert^2\leq \lVert x \rVert_+^2\leq \frac{\kappa}{\kappa_+}\lVert x \rVert^2 \qquad \forall\ x\in E.\]
\end{itemize}
The regularity constant $\kappa_E$ of the space
$(E,\lVert\cdot\rVert)$ is defined as:
\[\kappa_E:=\inf\{\kappa\geq 1:\ (E,\lVert\cdot\rVert)\text{ is }\kappa\text{-regular}\}.\]
From now on, assume that $(E,\lVert\cdot\rVert)$ has regularity
constant $\kappa_E$.

We define an $n$-dimensional martingale difference sequence
$\xi_1,\ldots,\xi_T$, that is, a sequence of $n$-dimensional random
vectors such that $\xi_{t-1}$ is $\sigma(\xi_t)$-measurable and
$\mathbb{E}_{\xi_t}\left\{\xi_t|\xi_{t-1}\right\} = 0$ for every $t$. In our context, the
following result on regular Euclidean spaces is of particular
interest.
\begin{thm}\cite{Nemirovski_Regular_Banach_Spaces}
\label{thm:Regular_Banach_Spaces} Choose $\chi\in (0,2]$ and reals
$\sigma_1,\ldots,\sigma_T>0$ such that:
\[
\mathbb{E}_{\xi_t}\left\{\exp\left\{\frac{\Vert
\xi_t\rVert^\chi}{\sigma_t^\chi }\right\}\arrowvert\
\xi_{t-1}\right\}\leq \exp\{ 1\}\qquad \forall\ t=1,\ldots,T.
\]
\begin{enumerate}
\item[(a) ] For all $c\geq 0$, we have:
\[
\mathbb{P}\left[\left\lVert\sum_{t=1}^T\xi_t\right\rVert>
c\sqrt{\kappa_E\sum_{t=1}^T\sigma_t^2} \right]\leq
C_\chi\exp\left(-c^\chi/C_\chi\right),
\]
where $C_\chi\geq 2$ is a properly chosen constant that solely
depends on $\chi$ and that is continuous in $\chi$. \item[(b) ] With
a properly chosen constant $c_\chi >0$ that solely depends on $\chi$
and that is continuous in $\chi$, we have:
\[\mathbb{E}_{\xi_{[T]}}\left\{\exp\left\{\left( \frac{ \lVert
\sum_{t=1}^T\xi_t\rVert}{c_\chi
\sqrt{\kappa_E\sum_{t=1}^T\sigma_t^2}}\right)^\chi\right\}\right\}\leq
\exp\left\{1\right\}.\]
\end{enumerate}
\end{thm}
\qed

For a proof of the first part of the above theorem, we refer to \cite{Nemirovski_Regular_Banach_Spaces}. Statement b) follows from a) by integration.

\section{Proofs}

\subsection{Proof of Theorem \ref{thm:convergence}}

The proof of Theorem \ref{thm:convergence} requires a result from
\cite{Juditsky08_solvingvariational}, which we reproduce below.

\begin{lem} \cite{Juditsky08_solvingvariational}
\label{thm:Lemma3.1inDetProx} Let $z\in Q^o$, $\gamma>0$, and $\eta,\zeta\in
E$. Consider the points
\begin{eqnarray*}
x&:=&\arg\min_{y\in Q}\left\lbrace \left\langle \gamma \eta-\omega'(z),y
\right\rangle +\omega(y) \right\rbrace \cr
z_+&:=&\arg\min_{y\in
Q}\left\lbrace \left\langle \gamma \zeta-\omega'(z),y \right\rangle
+\omega(y) \right\rbrace.
\end{eqnarray*}
Then,
\begin{eqnarray*}
\left\langle \gamma \zeta, x-y\right\rangle\leq
V_z(y)-V_{z_+}(y)+\frac{\gamma^2}{2}\left\|\eta-\zeta\right\|_{\ast}^2-\frac{1}{2}
\left\|x-z\right\|^2\qquad  \forall\ y\in Q.
\end{eqnarray*}
\end{lem}
\qed

Let us show Theorem \ref{thm:convergence}:
\\
\proof We can represent the random elements $F(z_{t-1})$ and $F(w_t)$
as follows:
\begin{eqnarray*}
\label{eq:splitting}
F(z_{t-1})&=&\hat{F}_{\xi_{2t-1}}(z_{t-1})-\sigma_{z_{t-1}}-\mu_{z_{t-1}}\qquad
\text{and} \qquad
F(w_t)=\hat{F}_{\xi_{2t}}(w_{t})-\sigma_{w_t}-\mu_{w_t},
\end{eqnarray*}
respectively. Let $u\in Q$. As $F$ is a monotone operator, we have:
\begin{eqnarray}
\label{eq:upper_bound_error} \sum_{t=1}^T\gamma_t\left\langle
F(u),z^T-u\right\rangle = \sum_{t=1}^T\gamma_t\left\langle
F(u),w_t-u\right\rangle &\leq& \sum_{t=1}^T\gamma_t\left\langle
F(w_t),w_t-u\right\rangle\cr &=&\sum_{t=1}^T\gamma_t\left\langle
\hat{F}_{\xi_{2t}}(w_{t})-\sigma_{w_t}-\mu_{w_t},w_t-u\right\rangle.
\end{eqnarray}
By Theorem \ref{thm:Lemma3.1inDetProx}, we obtain:
\begin{eqnarray*}
\sum_{t=1}^T\gamma_t\left\langle
\hat{F}_{\xi_{2t}}(w_{t}),w_t-u\right\rangle \leq
V_{z^\omega}(u)+\sum_{t=1}^T\left(\frac{\gamma_t^2}
{2}\left\|\hat{F}_{\xi_{2t}}(w_{t})-\hat{F}_{\xi_{2t-1}}(z_{t-1})\right\|_\ast^2-\frac{1}{2}\left\|w_{t}-z_{t-1}\right\|^2\right).
\end{eqnarray*}
Furthermore,
\begin{eqnarray*}
\left\|\hat{F}_{\xi_{2t}}(w_{t})-\hat{F}_{\xi_{2t-1}}(z_{t-1})\right\|_\ast^2&=&
\left\|F(w_t) +\sigma_{w_t} +\mu_{w_t}-F(z_{t-1})
-\sigma_{z_{t-1}}-\mu_{z_{t-1}}\right\|_\ast^2\cr &\leq&
2\left(\left\|F(w_t)-F(z_{t-1})\right\|_\ast^2 +
\left\|\sigma_{w_t}-\sigma_{z_{t-1}}+\mu_{w_t}-\mu_{z_{t-1}}\right\|_\ast^2\right)\cr
&\leq& 2\left(L^2 \left\|w_t-z_{t-1}\right\|^2+
\left\|\sigma_{w_t}-\sigma_{z_{t-1}} +
\mu_{w_t}-\mu_{z_{t-1}}\right\|_\ast^2\right),
\end{eqnarray*}
where the concluding inequality is due to the Lipschitz continuity of
$F$. Observe that $\gamma_t^2L^2\leq \frac{1}{2}$ because of the
step-size choice (\ref{eq:stepsize_bound}). Thus,
\begin{eqnarray}
\label{eq:intermediate_bound1} \sum_{t=1}^T\gamma_t\left\langle
\hat{F}_{\xi_{2t}}(w_{t}),w_t-u\right\rangle&\leq& V_{z^\omega}(u)
+ \sum_{t=1}^T\gamma_t^2
\left\|\sigma_{w_t}-\sigma_{z_{t-1}}+\mu_{w_t}
-\mu_{z_{t-1}}\right\|_\ast^2\cr &\leq&
\frac{\Omega^2}{2}+2\sum_{t=1}^T\gamma_t^2\left(\lVert
\sigma_{w_t}-\sigma_{z_{t-1}}\rVert_\ast^2+\lVert \mu_{w_t}+
\mu_{z_{t-1}}\rVert_\ast^2\right).
\end{eqnarray}
Additionally, the following inequality holds:
\begin{eqnarray*}
\label{eq:intermediate_bound2} \sum_{t=1}^T\gamma_t \left\langle
\sigma_{w_t}+\mu_{w_t},u-w_t\right\rangle &=& \sum_{t=1}^T\gamma_t
\left(\left\langle \sigma_{w_t},u-z^\omega\right\rangle+\left\langle
\sigma_{w_t},z^\omega-w_t\right\rangle + \left\langle
\mu_{w_t},u-w_t\right\rangle\right)\cr &\leq&  \Omega \left\|
\sum_{t=1}^T\gamma_t\sigma_{w_t}\right\|_\ast +
\sum_{t=1}^T\gamma_t\left( \left\langle
\sigma_{w_t},z^\omega-w_t\right\rangle + D\left\|
\mu_{w_t}\right\|_*\right).
\end{eqnarray*}
We combine the above inequality with (\ref{eq:upper_bound_error}) and (\ref{eq:intermediate_bound1}):
\begin{eqnarray*}
\sum_{t=1}^T\gamma_t\left\langle F(u),z^T-u\right\rangle &\leq&
\frac{\Omega^2}{2}+2\sum_{t=1}^T\gamma_t^2 \left(\lVert
\sigma_{w_t}-\sigma_{z_{t-1}}\rVert_\ast^2+\lVert \mu_{w_t}+
\mu_{z_{t-1}}\rVert_\ast^2\right)\cr &&+ \Omega \left\|
\sum_{t=1}^T\gamma_t\sigma_{w_t}\right\|_\ast +
\sum_{t=1}^T\gamma_t\left( \left\langle
\sigma_{w_t},z^\omega-w_t\right\rangle + D\left\|
\mu_{w_t}\right\|_*\right).
\end{eqnarray*}
Maximizing the left-hand side of the above inequality with respect to
$u\in Q$ and taking expectations on both sides, we obtain:
\begin{eqnarray*}
\sum_{t=1}^T\gamma_t\bE_{\xi_{[2T]}}\left\{\epsilon(z^T)\right\}
&\leq&  \frac{\Omega^2}{2} + 2\sum_{t=1}^T\gamma_t^2 \bE_{\xi_{[2T]}}
\left\{\lVert \sigma_{w_t}-\sigma_{z_{t-1}}\rVert_\ast^2+\lVert
\mu_{w_t}+ \mu_{z_{t-1}}\rVert_\ast^2\right\}\cr &&+ \Omega
\bE_{\xi_{[2T]}} \left\|
\sum_{t=1}^T\gamma_t\sigma_{w_t}\right\|_\ast + \sum_{t=1}^T\gamma_t
\bE_{\xi_{[2T]}} \left\{ \left\langle
\sigma_{w_t},z^\omega-w_t\right\rangle + D\left\|
\mu_{w_t}\right\|_*\right\}.
\end{eqnarray*}
It remains to observe that $\bE_{\xi_{[2T]}} \left\{ \left\langle
\sigma_{w_t},z^\omega-w_t\right\rangle\right\}=0$, as $\sigma_{w_t}$
is a martingale difference.

Assume that $F$ is associated with the saddle-point problem
(\ref{eq:saddle-point}), that is:
\[F(x,y)=\left(\frac{\partial
\phi(x,y)}{\partial x}; -\frac{\partial \phi(x,y)}{\partial
y}\right).\] Recall that $t\in\{1,\ldots,T\}$. Let $w_t=(x_t,y_t)\in
Q:=Q_1\times Q_2$, $u=(u_x,u_y)\in Q_1\times Q_2$, and
$\lambda_t:=\gamma_t/\sum_{i=1}^T\gamma_i$. As the function $\phi$ is
convex in the first and concave in the second argument, we have:
\begin{eqnarray*}
\sum_{t=1}^T\lambda_t\left\langle F(w_t),w_t-u\right\rangle &\geq&
\sum_{t=1}^T\lambda_t\left(\phi(x_t,y_t)-\phi(u_x,y_t)+\phi(x_t,u_y)
-\phi(x_t,y_t)\right)\cr &=&
\sum_{t=1}^T\lambda_t\left(\phi(x_t,u_y)-\phi(u_x,y_t)\right)\cr
&\geq& \phi\left(\sum_{t=1}^T\lambda_t
x_t,u_y\right)-\phi\left(u_x,\sum_{t=1}^T\lambda_t y_t\right).
\end{eqnarray*}
It remains to apply the same arguments as above in order to complete the proof.
\qed

\subsection{Proof of Proposition \ref{prop:matrix_exp_approximation}}

For $V\in \cS_n$, we define
\begin{eqnarray}
 \label{def:g_appendix}
g(V):=h(V)+c,\qquad \text{where } h(V):=\mathcal{A}^*\left(\frac{\exp\{V\}}{\text{Tr}(\exp\{V\})}\right).
\end{eqnarray}
Recall that $\xi^1,\ldots,\xi^N$ are independent
$\mathcal{N}(0,I_n)$-distributed random vectors and $\xi=(\xi^1,\ldots,\xi^N)$. Let  $g_\xi(V)$ be defined as in (\ref{estimate}), that is:
\begin{eqnarray}
 \label{def:g_xi_appendix}
 g_\xi(V):=h_\xi(V)+c,\qquad\text{where }h_\xi(V):=\mathcal{A}^*\left(\frac{G_\xi(V)}{\theta_\xi(V)}\right)
\end{eqnarray}
with
\begin{eqnarray}
\label{def:g_xi_appendix_2}
G_\xi(V):=\frac{\sum_{i=1}^N \left(\exp\left\{V/2\right\}\xi^i\right)\left(\exp\left\{V/2\right\}\xi^i\right)^T}{N}\qquad \text{and}\qquad \theta_\xi(V):=\frac{\sum_{i=1}^N[\xi^i]^T\exp\left\{V\right\}\xi^i}{N}.
\end{eqnarray}
We start with the observation:
\begin{ass}
As the standard multivariate normal distribution $\mathcal{N}(0,I_n)$ is orthogonal invariant, and as both
$h(V)$ and $h_\xi(V)$ are invariant under positive scaling, we can assume, without loss of
generality, that $\exp\{V\}$ is diagonal (with positive diagonal entries) and
of trace $1$. This assumption shall hold for the rest of this
section.
\end{ass}

For $i=1,\ldots,N$, we set:
\begin{eqnarray*}
D_{\xi^i}:=
\left(\exp\{V/2\}\xi^i\right)\left(\exp\{V/2\}\xi^i\right)^T-\exp\{V\},\qquad
d_{\xi^i}&:=&\cA^*D_{\xi^i},\qquad
d_{\xi}:=\frac{1}{N}\sum_{i=1}^Nd_{\xi^i},
\end{eqnarray*}
and:
\begin{eqnarray*}
f_{\xi^i}&:=&[\xi^i]^T\exp\{V\} \xi^i-1, \qquad
f_{\xi}:=\frac{1}{N}\sum_{i=1}^Nf_{\xi^i}.
\end{eqnarray*}

\begin{lem}
\label{lem:auxiliary_matrix_exp_bounds}
For an appropriately chosen constant $c>0$, we have for any $i=1,\ldots,N$.:
\begin{enumerate}
 \item[a) ]
$\bE_{\xi^i} D_{\xi^i}=\bE_{\xi^i} d_{\xi^i}=\bE_{\xi}
d_{\xi}= \bE_{\xi^i} f_{\xi^i}=\bE_{\xi} f_{\xi}=0$.
 \item[b) ]
$\bE_{\xi^i}\left\{ \exp\left\{\frac{\lVert
D_{\xi^i}\rVert_{Y}}{c}\right\}\right\} \leq \exp\{1\}$.
 \item[c) ]
$\bE_{\xi^i}\left\{ \exp\left\{ \frac{\lVert
d_{\xi^i}\rVert_{x,*}}{c\mathcal{L}}\right\}\right\} \leq
\exp\{1\}$.
 \item[d) ]
$\bE_{\xi} \left\{ \exp\left\{\frac{\sqrt{N} \lVert
d_{\xi}\rVert_{x,*}}{c\mathcal{L}\sqrt{\kappa}}\right\}\right\}
\leq \exp\{1\}$.
 \item[e) ]
$\bE_{\xi^i}\left\{ \exp\left\{\frac{\lvert
f_{\xi^i}\rvert}{c}\right\}\right\} \leq \exp\{1\}$.
 \item[f) ]
$\bE_{\xi}\left\{ \exp\left\{\frac{\sqrt{N} \lvert
f_{\xi}\rvert}{c}\right\}\right\} \leq \exp\{1\}$.
\end{enumerate}
\end{lem}

\proof
Assume that $diag(\exp\{V\})=(v_1,\ldots,v_n)^T$, where $v_i\geq 0$ for any $i=1,\ldots,n$ and $\sum_{i=1}^nv_i=1$.
\begin{enumerate}
 \item[a) ] Let $i\in\left\{1,\ldots,n\right\}$. We have:
\[\bE_{\xi^i}D_{\xi^i}=\exp\{V/2\}\bE_{\xi^i}\left\{\xi^i[\xi^i]^T\right\}\exp\{V/2\}-\exp\{V\}=0,\]
where the concluding equality holds as
$\xi^i\sim\mathcal{N}(0,I_n)$. Moreover,
\[
\bE_{\xi^i}\left\{[\xi^i]^T\exp\{V\}\xi^i\right\}=\sum_{k=1}^n v_k=1,
\]
which proves $\bE_{\xi^i} f_{\xi^i}=0$. The remaining equalities
follow immediately.
\item[b) ] Assume that the $n$-dimensional random vector $\zeta$ is $\mathcal{N}(0,I_n)$-distributed. Then,
\begin{eqnarray}
\label{eq:norm_of_H*eta_}
 \left\| \left( \exp\{V/2\}\zeta\right)\left( \exp\{V/2\}\zeta\right)^T\right\|_{Y} =\sum_{i=1}^nv_i\zeta_i^2.
\end{eqnarray}
For any $0< c_1< \frac{1-\exp\{-2\}}{2} \left(<\frac{1}{2}\right)$, it holds that:
\begin{eqnarray}
\label{eq:exp_bound_exp_on_norm_of_H_eta_} \bE_{\zeta} \exp\left\{c_1
\left\|\left(\exp\{V/2\}\zeta\right)\left(\exp\{V/2\}\zeta\right)^T\right\|_Y\right\}
&=&\prod_{i=1}^n\bE_{\zeta}\exp\left\{c_1 v_i\zeta_i^2\right\}\cr
&=&\prod_{i=1}^n\left(1-2c_1 v_i\right)^{-1/2}\cr
&=&\exp\left\{-\frac{1}{2}\sum_{i=1}^n\ln\left(1-2c_1
v_i\right)\right\}\cr &\leq&
\exp\left\{-\frac{1}{2}\ln\left(1-2c_1\right)\right\}\cr &\leq&\exp\{ 1\},
\end{eqnarray}
where the first inequality holds as the maximum of
$-\sum_{i=1}^n\ln\left(1-2c_1 v_i\right)$ over the probability
simplex is attained at an extreme point (note that the function
$v\mapsto -\sum_{i=1}^n\ln\left(1-2c_1 v_i\right)$ is separable and
each of its components is convex). We obtain:
\begin{eqnarray*}
\exp\{ 1\}\geq c_1
\bE_{\zeta}\left\|\left(\exp\{V/2\}\zeta\right)\left(\exp\{V/2\}\zeta\right)^T\right\|_Y.
\end{eqnarray*}
Using Jensen's inequality, it follows that:
\begin{eqnarray*}
\bE_{\zeta}\exp\left\{\frac{c_1}{2\exp\{ 1\}}\left\|\chi\chi^T
-\bE_{\zeta}\chi\chi^T\right\|_Y
\right\}\leq \exp\{ 1\}, \qquad \chi:=\exp\{V/2\}\zeta.
\end{eqnarray*}
We conclude by observing that
$\left(\exp\{V/2\}\zeta\right)\left(\exp\{V/2\}\zeta\right)^T
-\bE_{\zeta}\left(\exp\{V/2\}\zeta\right)\left(\exp\{V/2\}\zeta\right)^T$ has
the same distribution than any of the matrices $D_{\xi^i}$,
$i=1,\ldots,N$. \item[c) ] Let $i=1,\ldots,N$ and $0<
c_2\leq\frac{c_1}{2\exp\{ 1\}}$. Due to a), we obtain:
\begin{eqnarray*}
\bE_{\xi^i}\exp\left\{\frac{c_2\left\|d_{\xi^i}\right\|_{x,*}}{\mathcal{L}}\right\}=
\bE_{\xi^i}\exp\left\{\frac{c_2\left\|\cA^*D_{\xi^i}\right\|_{x,*}}{\mathcal{L}}\right\}\leq
\bE_{\xi^i}\exp\left\{\frac{c_2\left\|D_{\xi^i}\right\|_Y\mathcal{L}}{\mathcal{L}}\right\}\leq\exp\{ 1\}.
\end{eqnarray*}
\item[d) ] According to Theorem \ref{thm:Regular_Banach_Spaces}
(note that we apply Theorem \ref{thm:Regular_Banach_Spaces} with
$\chi=1$ and $\sigma_i=\mathcal{L}/c_2$ for any $i=1,\ldots,N$),
there exists $c_3>0$ such that:
\begin{eqnarray*}
\bE_{\xi}\exp\left\{\frac{c_2\sqrt{N}\left\|
d_{\xi}\right\|_{x,*}}{c_3\mathcal{L}\sqrt{\kappa}}\right\}\leq
\exp\{ 1\}.
\end{eqnarray*}
\item[e) ] Here, the $n$-dimensional random vector $\zeta$ is
$\mathcal{N}(0,I_n)$-distributed. Due to
(\ref{eq:norm_of_H*eta_}) and
(\ref{eq:exp_bound_exp_on_norm_of_H_eta_}), we obtain:
\begin{eqnarray*}
\bE_{\zeta}\exp\left\{\frac{c_1}{2}\left|\zeta^T \exp\{V\}\zeta-1 \right|
\right\}\leq\exp\left\{\frac{1}{2}\right\}\left( \bE_{\zeta}\exp\left\{c_1
\sum_{i=1}^nv_i\zeta_i\right\}\right)^{1/2}\leq \exp\{ 1\},
\end{eqnarray*}
where $0< c_1< \frac{1-\exp\{-2\}}{2}$. It remains to note that $\zeta^T
\exp\{V\}\zeta-1$ has the same distribution as any of the random variables
$f_{\xi^i}$, $i=1,\ldots,N$. \item[f) ] We observe that the space
$(\mathbb{R},\left\|\cdot\right\|_1)$ has a regularity constant of
$1$. Due to Theorem \ref{thm:Regular_Banach_Spaces} (we apply Theorem
\ref{thm:Regular_Banach_Spaces} with $\chi=1$ and $\sigma_i=2/c_1$
for any $i=1,\ldots,N$), there exists a constant $c_4>0$ such that:
\begin{eqnarray*}
\bE_{\xi}\exp\left\{\frac{c_1 \sqrt{N}\left| f_{\xi} \right|}{2c_4}
\right\}= \bE_{\xi}\exp\left\{\frac{c_1\left| \sum_{i=1}^N f_{\xi^i}
\right|}{2c_4\sqrt{N}} \right\}\leq\exp\{ 1\}.
\end{eqnarray*}
\end{enumerate}
It remains to choose $c\geq\max\left\lbrace 2\exp\{ 1\}/c_1,c_3/c_2,2c_4/c_1\right\rbrace$.
\qed

We are ready to prove Proposition \ref{prop:matrix_exp_approximation}:
\\
\proof Let $V\in\mathcal{S}_n$ and recall definitions (\ref{def:g_appendix})-(\ref{def:g_xi_appendix_2}). Consider the random element:
\[
\beta_{\xi}:=d_{\xi}-f_{\xi}h(V).
\]
Lemma \ref{lem:auxiliary_matrix_exp_bounds} implies
$\bE_{\xi}\beta_{\xi}=0$. Moreover, by the same lemma, there exists a
constant $c_1>0$ such that:
\begin{eqnarray}
\label{eq:bound_exp_beta} \bE_{\xi}
\exp\left\{\frac{\sqrt{N}\left\|\beta_{\xi}\right\|_{x,*}}{2c_1\mathcal{L}\sqrt{\kappa}}\right\}
&\leq& \bE_{\xi}\exp\left\{\frac{\sqrt{N}
\left\|d_{\xi}\right\|_{x,*}}{2c_1\mathcal{L}\sqrt{\kappa}}\right\}\exp\left\{\frac{\sqrt{N}
\left\|f_{\xi}h(V)\right\|_{x,*}}{2c_1\mathcal{L}\sqrt{\kappa}}\right\}\cr
&\leq& \left[\bE_{\xi}
\exp\left\{\frac{\sqrt{N}\left\|d_{\xi}\right\|_{x,*}}{c_1\mathcal{L}\sqrt{\kappa}}\right\}\right]^{1/2}
\left[\bE_{\xi}\exp\left\{\frac{\sqrt{N}\left\|f_{\xi}h(V)\right\|_{x,*}}{c_1\mathcal{L}\sqrt{\kappa}}
\right\}\right]^{1/2}\cr
&\leq&\exp\left\{\frac{1}{2}\right\}\left[\bE_{\xi}
\exp\left\{\frac{\sqrt{N}\left|f_{\xi}\right|}{c_1\sqrt{\kappa}}\right\}\right]^{1/2}\cr
&\leq&\exp\{ 1\},
\end{eqnarray}
where we use H\"older's inequality and the facts
$\left\|h(V)\right\|_{x,*}\leq\mathcal{L}$ and $\kappa\geq
1$. Let
\[
\gamma_{\xi}:=h_{\xi}(V)-\beta_{\xi}-h(V).
\]
As $\theta_{\xi}(V)=f_{\xi}+1$ and
$\cA^*\left(G_{\xi}(V)\right)=d_{\xi}+h(V)$, it holds that:
\[
h_{\xi}(V)=\frac{\cA^*\left(G_{\xi}(V) \right)}{\theta_{\xi}(V)}=\frac{d_{\xi}+h(V)}{f_{\xi}+1}.
\]
We obtain:
\[
\gamma_{\xi}=
\frac{d_{\xi}+h(V)}{f_{\xi}+1}-h(V)-d_{\xi}+f_{\xi}h(V)=\frac{f^2_{\xi}h(V)-d_{\xi}f_{\xi}}{1+f_{\xi}}.
\]
Consider the sets:
\[
\Pi:=\left\{\xi:\ \left|f_{\xi}\right|\leq
1/2\right\}\qquad\text{and}\qquad \hat{\Pi}:=\R^{n\times N}\setminus \Pi.
\]
When $\xi\in \Pi$, we have:
\begin{eqnarray}
\label{eq:bound_elements_of_Pi}
\left\|\gamma_{\xi}\right\|_{x,*}&=&\frac{\left\|
f^2_{\xi}h(V)-d_{\xi}f_{\xi}\right\|_{x,*}}{\left|1+f_{\xi}\right|}\leq
2\left\| f^2_{\xi}h(V)-d_{\xi}f_{\xi}\right\|_{x,*}\cr &\leq&
2\left(\left| f_{\xi}\right|^2\mathcal{L}+\left|
f_{\xi}\right|\left\|d_{\xi}\right\|_{x,*}\right)\leq\left|
f_{\xi}\right|\mathcal{L}+\left\|d_{\xi}\right\|_{x,*}.
\end{eqnarray}
More generally, it holds that:
\begin{eqnarray}
\label{eq:bound_elements_of_Pi_hat}
\left\|\gamma_{\xi}\right\|_{x,*}&\leq&
\left\|h_{\xi}(V)\right\|_{x,*}+\left\|\beta_{\xi}\right\|_{x,*}+\left\|h(V)\right\|_{x,*}\cr
&\leq&\left\|h_{\xi}(V)\right\|_{x,*} +
\left\|d_{\xi}\right\|_{x,*}+\left\|f_{\xi}h(V)\right\|_{x,*}+\left\|h(V)\right\|_{x,*}\cr
&=&
\left\|h_{\xi}(V)\right\|_{x,*}+\left\|d_{\xi}\right\|_{x,*}+
\left|f_{\xi}\right|\left\|h(V)\right\|_{x,*}+\left\|h(V)\right\|_{x,*}\cr
&\leq&
\left(2+\left|f_{\xi}\right|\right)\mathcal{L}+\left\|d_{\xi}\right\|_{x,*},
\end{eqnarray}
which is due to $\left\|h(V)\right\|_{x,*}\leq\mathcal{L}$ and to
the following inequality:
\begin{eqnarray*}
\left\|h_{\xi}(V) \right\|_{x,*}&\leq& \mathcal{L}
\frac{\left\|\sum_{i=1}^N\left(\exp\{V/2\}\xi^i
\right)\left(\exp\{V/2\}\xi^i
\right)^T\right\|_Y}{\sum_{i=1}^N[\xi^i]^T \exp\{V\}
\xi^i}=\mathcal{L}.
\end{eqnarray*}
Furthermore, denoting by $\mathbb{P}$ the probability measure of the
random matrix $\xi$:
\begin{eqnarray}
\label{eq:prob_bound_on_Pi_hat} \mathbb{P}\left[\hat{\Pi}\right]&=&
\mathbb{P}\left[\left|f_{\xi}\right|>1/2\right] =
\mathbb{P}\left[\exp\left\{\frac{\sqrt{N}\left|f_{\xi}\right|}{c_1}\right\}
>\exp\left\{\frac{\sqrt{N}}{2c_1}\right\}\right]\cr
&<&\exp\left\{-\frac{\sqrt{N}}{2c_1}\right\}\bE_{\xi}
\exp\left\{\frac{\sqrt{N}\left|f_{\xi}\right|}{c_1}\right\}\cr
&\leq&\exp\left\{-\frac{\sqrt{N}}{2c_1}\right\}\exp\left\{ 1\right\},
\end{eqnarray}
where the inequalities follow from Markov's inequality and from
statement e) of Lemma \ref{lem:auxiliary_matrix_exp_bounds},
respectively. Choose $c_2\geq 4c_1$ and observe:
\begin{eqnarray*}
\bE_{\xi}\exp\left\{\frac{\sqrt{N}\left\|\gamma_{\xi}\right\|_{x,*}}
{3c_2\mathcal{L}\sqrt{\kappa}}\right\}&=& \int_{\Pi}\exp\left\{
\frac{\sqrt{N}\left\|\gamma_{\xi}\right\|_{x,*}}
{3c_2\mathcal{L}\sqrt{\kappa}}\right\}d\mathbb{P}(\xi) +
\int_{\hat{\Pi}}\exp\left\{
\frac{\sqrt{N}\left\|\gamma_{\xi}\right\|_{x,*}}
{3c_2\mathcal{L}\sqrt{\kappa}}\right\}d\mathbb{P}(\xi).
\end{eqnarray*}
By (\ref{eq:bound_elements_of_Pi}), Cauchy-Schwarz inequality, and
Lemma \ref{lem:auxiliary_matrix_exp_bounds}, we obtain:
\begin{eqnarray*}
\int_{\Pi}\exp\left\{
\frac{\sqrt{N}\left\|\gamma_{\xi}\right\|_{x,*}}
{3c_2\mathcal{L}\sqrt{\kappa}}\right\}d\mathbb{P}(\xi)&\leq&
\int_{\Pi}\exp\left\{ \frac{\sqrt{N}\left(\left|
f_{\xi}\right|\mathcal{L}+ \left\|d_{\xi}\right\|_{x,*}\right)}
{3c_2\mathcal{L}\sqrt{\kappa}}\right\}d\mathbb{P}(\xi)\cr
&\leq&1\cdot \left[\bE_{\xi} \exp\left\{\frac{\sqrt{N}\left|
f_{\xi}\right|}{c_2\sqrt{\kappa}}\right\}\right]^{1/3}
\left[\bE_{\xi}\exp\left\{
\frac{\sqrt{N}\left\|d_{\xi}\right\|_{x,*}}
{c_2\mathcal{L}\sqrt{\kappa}}\right\}\right]^{1/3}\cr &\leq&
\exp\left\{\frac{1}{3}\right\}\exp\left\{\frac{1}{3}\right\}
=\exp\left\{\frac{2}{3}\right\}.
\end{eqnarray*}
Additionally, by (\ref{eq:bound_elements_of_Pi_hat}), we have:
\begin{eqnarray*}
\int_{\hat{\Pi}}\exp\left\{
\frac{\sqrt{N}\left\|\gamma_{\xi}\right\|_{x,*}}
{3c_2\mathcal{L}\sqrt{\kappa}}\right\}d\mathbb{P}(\xi)&\leq&
\int_{\hat{\Pi}}\exp\left\{
\frac{\sqrt{N}\left(2\mathcal{L}+\left|f_{\xi}\right|\mathcal{L} +
\left\|d_{\xi}\right\|_{x,*}\right)}
{3c_2\mathcal{L}\sqrt{\kappa}}\right\}d\mathbb{P}(\xi).
\end{eqnarray*}
Let
\[
A:= \int_{\hat{\Pi}}\exp\left\{
\frac{\sqrt{N}\left(\left|f_{\xi}\right|\mathcal{L} +
\left\|d_{\xi}\right\|_{x,*}\right)}
{3c_2\mathcal{L}\sqrt{\kappa}}\right\}d\mathbb{P}(\xi).
\]
Cauchy-Schwarz inequality, bound (\ref{eq:prob_bound_on_Pi_hat}), and
Lemma \ref{lem:auxiliary_matrix_exp_bounds} imply:
\begin{eqnarray*}
A&\leq& \left[\int_{\hat{\Pi}}1d\mathbb{P}(\xi)\right]^{1/3} \left[
\int_{\hat{\Pi}}\exp\left\{
\frac{\sqrt{N}\left|f_{\xi}\right|}{c_2\sqrt{\kappa}}\right\}
d\mathbb{P}(\xi)\right]^{1/3}\left[\int_{\hat{\Pi}}\exp\left\{
\frac{\sqrt{N} \left\|d_{\xi}\right\|_{x,*}}
{c_2\mathcal{L}\sqrt{\kappa}}\right\}d\mathbb{P}(\xi)\right]^{1/3}\cr
&\leq&\left[\mathbb{P}\left[\hat{\Pi}\right]\right]^{1/3}
\left[\bE_{\xi} \exp\left\{
\frac{\sqrt{N}\left|f_{\xi}\right|}{c_2\sqrt{\kappa}}\right\}
\right]^{1/3} \left[\bE_{\xi}\exp\left\{
\frac{\sqrt{N}\left\|d_{\xi}\right\|_{x,*}}{c_2\mathcal{L}
\sqrt{\kappa}}\right\}\right]^{1/3}\cr
&\leq&\left[\mathbb{P}\left[\hat{\Pi}\right]\right]^{1/3}
\left[\bE_{\xi}\exp\left\{
\frac{\sqrt{N}\left|f_{\xi}\right|}{c_2}\right\}
\right]^{1/3}\left[\bE_{\xi}\exp\left\{
\frac{\sqrt{N}\left\|d_{\xi}\right\|_{x,*}}{c_2\mathcal{L}\sqrt{\kappa}}\right\}\right]^{1/3}\cr
&\leq&
\exp\left\{-\frac{\sqrt{N}}{6c_1}\right\}\exp\left\{\frac{2}{3}\right\}.
\end{eqnarray*}
As $c_2\geq 4c_1$, we obtain:
\begin{eqnarray}
\label{eq:bound_exp_gamma} \bE_{\xi}\exp\left\{ \frac{\sqrt{N}
\left\|\gamma_{\xi}\right\|_{x,*}}{3c_2\mathcal{L}\sqrt{\kappa}}\right\}
&\leq&\exp\left\{\frac{2}{3}\right\}+\exp\left\{\frac{2\sqrt{N}}
{3c_2\sqrt{\kappa}}\right\}\exp\left\{-\frac{\sqrt{N}}{6c_1}\right\}\exp\left\{\frac{2}{3}\right\}\cr
&\leq&\exp\left\{\frac{2}{3}\right\}+\exp\left\{\frac{2\sqrt{N}}
{3c_2}-\frac{\sqrt{N}}{6c_1}\right\}\exp\left\{\frac{2}{3}\right\}\leq2\exp\left\{\frac{2}{3}\right\}.
\end{eqnarray}
Thus,
\begin{eqnarray}\label{eq:exp_g_hat-g}
\bE_{\xi}\exp\left\{\frac{\sqrt{N}\left\|h_{\xi}(V)-h(V)
\right\|_{x,*}}{6c_2\mathcal{L}\sqrt{\kappa}}\right\}&=&
\bE_{\xi}\exp\left\{\frac{\sqrt{N}\left\|\beta_{\xi}+
\gamma_{\xi}\right\|_{x,*}}{6c_2\mathcal{L}\sqrt{\kappa}}\right\}\cr
&\leq&\left[\bE_{\xi}\exp\left\{\frac{\sqrt{N}\left\|\beta_{\xi}\right\|_{x,*}}
{3c_2\mathcal{L}\sqrt{\kappa}}\right\}\right]^{1/2}\left[\bE_{\xi}
\exp\left\{\frac{\sqrt{N}\left\|\gamma_{\xi}\right\|_{x,*}}{3c_2\mathcal{L}\sqrt{\kappa}}\right\}\right]^{1/2}\cr
&\leq&\left[\bE_{\xi}\exp\left\{\frac{\sqrt{N}\left\|\beta_{\xi}\right\|_{x,*}}
{c_1\mathcal{L}\sqrt{\kappa}}\right\}\right]^{1/2}
\left[\bE_{\xi}\exp\left\{\frac{\sqrt{N}\left\|\gamma_{\xi}\right\|_{x,*}}
{3c_2\mathcal{L}\sqrt{\kappa}}\right\}\right]^{1/2}\cr
&\leq&\sqrt{2}\exp\left\{\frac{5}{6}\right\},
\end{eqnarray}
where the inequalities are due to H\"older's inequality, the fact
that $c_2\geq c_1$, bound (\ref{eq:bound_exp_beta}), and inequality
(\ref{eq:bound_exp_gamma}), respectively.

It remains to find an appropriate bound on the norm of
$h(V)-\bE_{\xi}h_{\xi}(V)$. Recall that
$\bE_{\xi}\beta_{\xi}=0$. Therefore,
\begin{eqnarray*}
\left\|\bE_{\xi}h_{\xi}(V)-h(V) \right\|_{x,*}&=&
\left\|\bE_{\xi}\gamma_{\xi}\right\|_{x,*}\leq \bE_{\xi}
\left\|\gamma_{\xi}\right\|_{x,*}\cr
&=&\int_{\Pi}\left\|\gamma_{\xi} \right\|_{x,*}d\mathbb{P}(\xi) +
\int_{\hat{\Pi}}\left\|\gamma_{\xi}\right\|_{x,*}d\mathbb{P}(\xi)\cr
&\leq&2 \int_{\Pi}\left| f_{\xi}\right|^2\mathcal{L}+\left|
f_{\xi}\right|\left\|d_{\xi}\right\|_{x,*}
d\mathbb{P}(\xi)+\int_{\hat{\Pi}}\left\|\gamma_{\xi}\right\|_{x,*}d\mathbb{P}(\xi),
\end{eqnarray*}
where the concluding inequality follows from
(\ref{eq:bound_elements_of_Pi}). As $2\exp(x)\geq x^2$ for any $x\geq
0$, we obtain by Lemma \ref{lem:auxiliary_matrix_exp_bounds}:
\begin{eqnarray*}
2\mathcal{L} \int_{\Pi}\left| f_{\xi}\right|^2d\mathbb{P}(\xi)&\leq&
2\mathcal{L}\bE_{\xi}\left|
f_{\xi}\right|^2\leq\frac{4c_1^2\mathcal{L}}{N}\bE_{\xi}\exp\left\{\frac{\sqrt{N}\left|
f_{\xi}\right|}{c_1}\right\}\leq\frac{4\exp\{ 1\}c_1^2\mathcal{L}}{N}.
\end{eqnarray*}
Furthermore, the same arguments imply:
\begin{eqnarray*}
2 \int_{\Pi}\left|
f_{\xi}\right|\left\|d_{\xi}\right\|_{x,*}d\mathbb{P}(\xi)&\leq&
\int_{\Pi} \sqrt{\kappa}\mathcal{L}\left|
f_{\xi}\right|^2d\mathbb{P}(\xi)+\int_{\Pi}\frac{\left\|d_{\xi}\right\|_{x,*}^2}
{\sqrt{\kappa}\mathcal{L}}d\mathbb{P}(\xi)\cr
&\leq&\bE_{\xi}\left\{\sqrt{\kappa}\mathcal{L}\left|
f_{\xi}\right|^2\right\}+\bE_{\xi}
\left\{\frac{\left\|d_{\xi}\right\|_{x,*}^2}{\mathcal{L}\sqrt{\kappa}}\right\}\cr
&\leq&\frac{2c_1^2\mathcal{L}\sqrt{\kappa}}{N}\left(\bE_{\xi}
\left\{\exp\left\{\frac{\sqrt{N}\left|f_{\xi}\right|}{c_1}\right\}\right\}
+ \bE_{\xi}\left\{\exp\left\{\frac{\sqrt{N}
\left\|d_{\xi}\right\|_{x,*}}{c_1\mathcal{L}\sqrt{\kappa}}\right\}\right\}\right)\cr
&\leq&\frac{4\exp\{ 1\}c_1^2\mathcal{L}\sqrt{\kappa}}{N}.
\end{eqnarray*}
Finally,
\begin{eqnarray*}
 \int_{\hat{\Pi}}\left\|\gamma_{\xi}\right\|_{x,*}d\mathbb{P}(\xi)
&\leq&
\left[\mathbb{P}\left[\hat{\Pi}\right]\right]^{1/2}\left[\bE_{\xi}
\left\|\gamma_{\xi}\right\|_{x,*}^2\right]^{1/2}\cr
&\leq&\exp\left\{
\frac{1}{2}-\frac{\sqrt{N}}{4c_1}\right\}\left[\bE_{\xi}
\left\|\gamma_{\xi}\right\|_{x,*}^2\right]^{1/2}\cr
&\leq&\exp\left\{
\frac{1}{2}-\frac{\sqrt{N}}{4c_1}\right\}\left[\frac{18c_2^2\mathcal{L}^2
\kappa}{N}\bE_{\xi}\exp\left\{\frac{\sqrt{N}
\left\|\gamma_{\xi}\right\|_{x,*}}{3c_2\mathcal{L}\sqrt{\kappa}}\right\}\right]^{1/2}\cr
&\leq&\frac{6\exp\left\{\frac{5}{6}\right\}c_2\mathcal{L}\sqrt{\kappa}}{\sqrt{N}}
\exp\left\{-\frac{\sqrt{N}}{4c_1}\right\}\cr
&\leq&\frac{6\exp\left\{\frac{5}{6}\right\}c_2\mathcal{L}\sqrt{\kappa}}{\sqrt{N}}
\exp\left\{-\frac{\sqrt{N}}{c_2}\right\}\cr
&\leq&\frac{6\exp\left\{\frac{5}{6}\right\}c_2^2\mathcal{L}\sqrt{\kappa}}{N},
\end{eqnarray*}
where the inequalities hold due to H\"older's inequality, the fact
that $\left\|\gamma_{\xi}\right\|_{x,*}^2$ is nonnegative
for any $\xi\in\mathbb{R}^{n\times N}$, bound (\ref{eq:prob_bound_on_Pi_hat}), the fact
that $2\exp(x)\geq x^2$ for any $x\geq 0$, inequality
(\ref{eq:bound_exp_gamma}), and the assumption $c_2\geq 4 c_1$,
respectively. We obtain:
\begin{eqnarray*}
\left\|\bE_{\xi}h_{\xi}(V)-h(V)\right\|_{x,*}&\leq&
\frac{4\exp\{ 1\}c_1^2\mathcal{L}}{N}+
\frac{4\exp\{ 1\}c_1^2\mathcal{L}\sqrt{\kappa}}{N}+\frac{6\exp\{ 1\}c_2^2\mathcal{L}\sqrt{\kappa}}{N}\cr
&\leq&
\frac{8\exp\{ 1\}c_1^2\mathcal{L}\sqrt{\kappa}}{N}+\frac{6\exp\{ 1\}c_2^2\mathcal{L}\sqrt{\kappa}}{N}\cr
&\leq&
\frac{\exp\{ 1\}c_2^2\mathcal{L}\sqrt{\kappa}}{2N}+\frac{6\exp\{ 1\}c_2^2\mathcal{L}\sqrt{\kappa}}{N}\cr
&=& \frac{13\exp\{ 1\}c_2^2\mathcal{L}\sqrt{\kappa}}{2N},
\end{eqnarray*}
which proves statement a). The above inequality ensures together with bound (\ref{eq:exp_g_hat-g}):
\begin{eqnarray*}
\bE_{\xi}\exp\left\{\frac{\sqrt{N}
\left\|h_{\xi}(V)-\bE_{\xi}h_{\xi}(V)\right\|_{x,*}}
{6c_2\mathcal{L}\sqrt{\kappa}}\right\}
&\leq&\sqrt{2}\exp\left\{\frac{5}{6}\right\}\exp\left\{\frac{13\exp\{ 1\}c_2}{12\sqrt{N}}\right\}\cr
&\leq& \sqrt{2}\exp\left\{\frac{5}{6}+\frac{13\exp\{ 1\}c_2}{12}\right\}.
\end{eqnarray*}
Let:
\begin{eqnarray*}
 c_3:=\frac{5}{6}+\frac{\ln(2)}{2}+\frac{13\exp\{ 1\}c_2}{12}.
\end{eqnarray*}
By Jensen's inequality, we conclude that:
\begin{eqnarray*}
\bE_{\xi}\exp\left\{\frac{\sqrt{N}
\left\|h_{\xi}(V)-\bE_{\xi}h_{\xi}(V)\right\|_{x,*}}
{6c_2c_3\mathcal{L}\sqrt{\kappa}}\right\} &\leq&\exp\{ 1\}.
\end{eqnarray*}
\qed

\end{document}